\documentclass[11pt, a4paper]{article}

\usepackage[english]{babel}
\usepackage{dutchcal}

\usepackage{amsmath}
\usepackage{amssymb}
\usepackage{amsthm}
\usepackage{mathabx}
\usepackage{tikz}
\usepackage{tikz-cd}
\usepackage{todonotes}
\usepackage{hyperref}
\usepackage{cleveref}
\usepackage{footmisc}
\usepackage{bbm}
\usepackage{afterpage}
\usepackage{mathrsfs}
\usepackage{sseq}
\usepackage{extarrows}
\usepackage{enumitem}
\usepackage{stmaryrd}
\usepackage{pdflscape}

\newcommand{\ab}{\mathrm{ab}}
\newcommand{\aff}{\mathrm{aff}}
\newcommand{\afp}{\mathrm{afp}}
\newcommand{\cn}{\mathrm{cn}}

\newcommand{\Ell}{\mathrm{Ell}}
\newcommand{\Temp}{\mathrm{Temp}}

\newcommand{\et}{\mathrm{\acute{e}t}}
\newcommand{\fin}{\mathrm{fin}}
\newcommand{\fgf}{\mathrm{fgf}}

\newcommand{\gl}{\mathrm{gl}}

\newcommand{\GM}{\mathrm{GM}}

\newcommand{\lex}{\mathrm{lex}}

\newcommand{\ori}{\mathrm{or}}
\newcommand{\res}{\mathrm{res}}

\newcommand{\BT}{{{\mathrm{BT}}}}

\DeclareMathOperator{\Spec}{Spec}

\newcommand{\Tors}{\mathrm{Tors}}

\DeclareMathOperator{\KU}{KU}

\DeclareMathOperator{\TMF}{TMF}

\DeclareMathOperator{\Ab}{Ab}
\DeclareMathOperator{\Aff}{Aff}

\DeclareMathOperator{\Cat}{Cat}
\DeclareMathOperator{\CAlg}{CAlg}

\DeclareMathOperator{\Fun}{Fun}
\DeclareMathOperator{\Glo}{Glo}

\DeclareMathOperator{\Mod}{Mod}

\newcommand{\Lie}{\mathrm{Lie}}

\newcommand{\OS}{\mathcal{OS}}
\DeclareMathOperator{\Pre}{Pre}
\DeclareMathOperator{\PreAb}{PreAb}
\newcommand{\PregeoAb}{\mathrm{PreAb}^{\mathrm{geo}}}
\newcommand{\geoAb}{\mathrm{Ab}^{\mathrm{geo}}}

\DeclareMathOperator{\PreBT}{PreBT}

\DeclareMathOperator{\QCoh}{QCoh}

\DeclareMathOperator{\Ran}{Ran}

\DeclareMathOperator{\Shv}{Shv}
\DeclareMathOperator{\Stk}{Stk}
\DeclareMathOperator{\Sp}{Sp}

\newcommand{\Spc}{\mathcal{S}}

\newcommand{\colim}{\mathrm{colim}\,}

\newcommand{\id}{\mathrm{id}}
\DeclareMathOperator{\Map}{Map}
\newcommand{\MAP}{\underline{\mathrm{Map}}}
\DeclareMathOperator{\Hom}{Hom}
\newcommand{\HOM}{\underline{\mathrm{Hom}}}
\newcommand{\op}{\mathrm{op}}

\newcommand{\A}{\mathsf{A}}
\newcommand{\B}{\mathbf{B}}
\renewcommand{\b}{\mathsf{B}}

\newcommand{\C}{\mathcal{C}}
\newcommand{\D}{\mathcal{D}}
\newcommand{\E}{\mathbf{E}}
\newcommand{\e}{\mathsf{E}}

\newcommand{\FF}{{\mathscr{F}}}
\newcommand{\G}{\mathbf{G}}
\newcommand{\GG}{\mathscr{G}}

\newcommand{\M}{\mathsf{M}}

\renewcommand{\O}{\mathcal{O}}

\renewcommand{\P}{\mathbf{P}}
\newcommand{\Q}{\mathbf{Q}}

\newcommand{\s}{\mathsf{S}}
\newcommand{\T}{\mathbf{T}}
\renewcommand{\t}{\mathsf{T}}
\newcommand{\x}{\mathsf{X}}
\newcommand{\xx}{\mathsf{X}}
\newcommand{\y}{\mathsf{Y}}
\newcommand{\Z}{\mathbf{Z}}

\newcommand{\Ga}{\Gamma}

\usepackage{mathtools}

\theoremstyle{theorem}\numberwithin{equation}{section}
\newtheorem{theorem}[equation]{Theorem}
\crefname{theorem}{{theorem}}{{theorems}}
\Crefname{theorem}{{Theorem}}{{Theorems}}
\newtheorem{theoremalph}{Theorem}

\crefname{theoremalph}{{theorem}}{{theorems}}
\Crefname{theoremalph}{{Theorem}}{{Theorems}}
\newtheorem{conjecturealph}[theoremalph]{Conjecture}

\Crefname{problem}{{Problem}}{{Problems}}
\newtheorem{prop}[equation]{Proposition}
\Crefname{prop}{{Proposition}}{{Propositions}}
\newtheorem{lemma}[equation]{Lemma}
\Crefname{lemma}{{Lemma}}{{Lemmas}}

\Crefname{cor}{{Corollary}}{{Corollaries}}

\Crefname{conjecture}{{Conjecture}}{{Conjectures}}

\theoremstyle{definition}\numberwithin{equation}{section}
\newtheorem{mydef}[equation]{Definition}
\Crefname{mydef}{{Definition}}{{Definitions}}

\Crefname{recall}{{Recollection}}{{Recollections}}

\Crefname{construction}{{Construction}}{{Constructions}}

\Crefname{ass}{{Assumptions}}{{Assumptions}}
\newtheorem{notation}[equation]{Notation}
\Crefname{notation}{{Notation}}{{Notations}}

\Crefname{situation}{{Situation}}{{Situations}}

\theoremstyle{remark}\numberwithin{equation}{section}
\newtheorem{example}[equation]{Example}
\Crefname{example}{{Example}}{{Examples}}

\Crefname{nonexample}{{NonExample}}{{NonExamples}}

\Crefname{claim}{{Claim}}{{Claims}}
\newtheorem{remark}[equation]{Remark}
\Crefname{remark}{{Remark}}{{Remarks}}

\Crefname{idea}{{Idea}}{{Ideas}}
\newtheorem{warn}[equation]{Warning}
\Crefname{warn}{{Warning}}{{Warnings}}

\Crefname{figure}{{Figure}}{{Figures}}
\Crefname{footnote}{{Footnote}}{{Footnotes.}}

\Crefname{part}{{\textsection}\!\!}{{\textsection}\!\!}
\Crefname{chapter}{{\textsection}\!\!}{{\textsection}\!\!}
\Crefname{section}{{\textsection}\!\!}{{\textsection}\!\!}
\Crefname{subsection}{{\textsection}\!\!}{{\textsection}\!\!}
\Crefname{appendix}{{\textsection}\!\!}{{\textsection}\!\!}

\begin{document}
\title{Comparing tempered and equivariant elliptic cohomology}
\author{Jack Morgan Davies\footnote{\url{davies@uni-wuppertal.de}}}
\maketitle

\begin{abstract}
In \cite{ec3} and \cite{davidandlennart}, Lurie and Gepner--Meier each define equivariant cohomology theories, namely \emph{tempered cohomology} and \emph{equivariant elliptic cohomology}, respectively, using derived algebraic geometry. We establish a natural equivalence between these theories on their overlap. Moreover, we emphasise the naturality and coherence of both these equivariant theories as well as our comparison. To demonstrate the use of this comparison, we show that the $G$-fixed points of equivariant topological modular forms is dualisable as a $\TMF$-module for all compact Lie groups $G$ that decompose as a product of a torus and a finite group. This is done by using our comparison to reduce to an argument of Gepner--Meier based on Lurie's tempered ambidexterity.
\end{abstract}

\setcounter{tocdepth}{3}
\tableofcontents

\newpage

\addcontentsline{toc}{section}{Introduction}
\section*{Introduction}
Equivariant cohomology theories have received much attention in recent years due to breakthroughs in computational techniques, theoretic frameworks, and applications, most famously demonstrated in \cite{hhr}. Similar to the classical nonequivariant situation, despite numerous tools and techniques now available to study stable equivariant stable homotopy theory, it is remarkably hard to produce genuinely new examples of equivariant cohomology theories other than variants of singular cohomology, $K$-theory, or geometric bordism.\\

In \cite{ec3}, Lurie constructs new families of equivariant cohomology theories which he coins \emph{tempered cohomology theories}. In more detail, given a \emph{preoriented $\P$-divisible group} $\G$ over an $\E_\infty$-ring $\A$, Lurie defines a functor $\A_\G\colon \OS_{\ab}^\op\to \CAlg_\A$ from the $\infty$-category of \emph{finite abelian orbispaces} to $\E_\infty$-$\A$-algebras. The functor $\A_\G$ is supposed to behave like an equivariant multiplicative cohomology theory defined for all finite groups simultaneously. For example, the assignment sending an $H$-space $X$ to the orbispace $X//H$ and then applying $\A_\G$ acts like an $H$-equivariant cohomology theory.\\

In \cite{davidandlennart}, Gepner--Meier use another technical set-up to construct equivariant cohomology theories from \emph{preoriented elliptic curves}, and more generally, from \emph{preoriented abelian sheaves}.\footnote{The first version of this article restricted its attention to the full subcategory of preoriented abelian sheaves spanned by those represented by abelian \emph{varieties}. The current generality in this article is built to further encompass geometric objects such as the multiplicative group scheme $\G_m$, which is closely tied to equivariant topological $K$-theory.} Given a preoriented abelian sheaf $\FF$ over a derived stack $\s$, they define a functor $\Ell_{\FF/\s}\colon \OS_\Lie\to \Shv(\s)$ from the $\infty$-category of \emph{compact Lie orbispaces} to sheaves on $\s$.

\subsection*{Comparison statements}

To compare these two definitions, we restrict our attention to \emph{geometric} preoriented abelian sheaves $\FF$, so those abelian sheaves whose associated torsion object $\FF[\P^\infty]$ is a $\P$-divisible group; see \Cref{torsionofsheavesdefintiion}. We write $\PregeoAb(\s)$ for the $\infty$-category of geometric preoriented abelian sheaves over a derived stack $\s$. The inclusion $\FF[\P^\infty]\to \FF$ then serves as the basis of our comparison between $\Ell_{\FF/\s}$ and $\Temp_{\FF[\P^\infty]/\s}$; the latter is our notation for the extension of Lurie's tempered cohomology from affine derived stacks $\Spec \A$ to general derived stacks $\s$ and extended from $\OS_\ab$ to orbispaces based on \emph{all} finite groups $\OS$. After showing that $\Temp_{[\P^\infty]}$ and $\Ell$ are functorial in the pair $\FF/\s$, see \Cref{functorialell,functorialtemp}, our comparison map is then defined as a natural transformation between these two equivariant cohomology theories, themselves defined as functors between Cartesian fibrations over the $\infty$-category of derived stacks $\Stk$.

\begin{theoremalph}\label{maintheoremgeo}
There is a natural equivalence of functors
\[\Temp_{[\P^\infty]}\xrightarrow{\simeq} \Ell\colon \int_{\Stk} \PregeoAb \to \int_{\Stk} \Fun^L(\OS,\Shv(-))\]
induced by the inclusion $\FF[\P^\infty]\to \FF$.
\end{theoremalph}

This first comparison concerns a ``geometric'' extension of tempered cohomology and equivariant elliptic cohomology. Conversely, we have an ``algebraic'' comparison statement between theories with values in algebras of quasi-coherent sheaves.

\begin{theoremalph}\label{maintheoremalg}
There is a natural equivalence of functors
\[\O_\Ell\xrightarrow{\simeq} \O_{\Temp_{[\P^\infty]}} \colon \int_{\Stk}\PregeoAb\to \int_{\Stk}\Fun^R(\OS^\op,\CAlg(\QCoh(-))\]
induced by the inclusion $\FF[\P^\infty]\to \FF$.
\end{theoremalph}

The global sections of $\O_{\Temp_{[\P^\infty]}}$ evaluated on a pair $\FF/\Spec \A$ is precisely Lurie's tempered cohomology construction $\A_{\FF[\P^\infty]}$ right Kan extended from $\OS_\ab$ to $\OS$, courtesy of the natural equivalence symmetric monoidal equivalence of $\infty$-categories $\QCoh(\Spec \A) \simeq \Mod_\A$.\\

As an example of the above theorems, consider the multiplicative group scheme $\G_m$ over $\KU$. \Cref{maintheoremgeo,maintheoremalg} then imply that the associated equivariant elliptic cohomology $\Ell_{\G_m/\KU}$ restricted to finite groups agrees with the tempered cohomology $\Temp_{\mu_{\P^\infty}/\KU}$ of the associated $\P$-divisible group $\mu_{\P^\infty} = \G_m[\P^\infty]$. Both \cite[\textsection4]{davidandlennart} and \cite[\textsection4.1]{ec3} independently justify that these constructions agree with the more classical equivariant complex $K$-theory of \cite{equivktheory}, however, the comparison afforded by the above theorems also ensures the \emph{naturality} of such a comparison. For instance, this implies that equivariant theories associated with \emph{real} topological $K$-theory, defined by $\G_m$ over $\Spec \KU / C_2$, are also equivalent, and that these equivalences in the real and complex case commute.\\

Other explicit examples come from elliptic curves over derived stacks. The initial case of this is the universal oriented elliptic curve $\e^\ori$ over the derived moduli stack of oriented elliptic curves $\M_\Ell^\ori$ of Lurie \cite[\textsection7]{ec2}. For a finite group $G$, one could call either of the $\E_\infty$-rings $\Ga(\Ell^{\B G}_{\e^\ori/\M_\Ell^\ori})$ or $\Ga(\Temp_{\e^\ori[\P^\infty]/\M_\Ell^\ori}^{\B G})$ the \emph{(genuine) $G$-fixed points of equivariant topological modular forms} $\TMF^{\B G}$. Applying global sections to \Cref{maintheoremalg}, we obtain a natural equivalence between these $\E_\infty$-rings, justifying the uniform notation $\TMF^{\B G}$. Again, the naturality of the above statements ensures the compatibility of this universal comparison with any equivariant elliptic cohomology theory. That includes those equivariant theories associated with equivariant topological modular forms with level structure or equivariant Tate $K$-theory; the former is explored in \cite{normsontmf} and the latter is constructed in \cite{globaltate}.\\

The proofs of both \Cref{maintheoremgeo} and \Cref{maintheoremalg} are formal consequences of the crucial case of restricting the functor $\infty$-categories from $\OS$ to $\Glo_\ab$, defined as the full subcategory of spaces spanned by those of the form $\B G$ for a finite abelian group $G$; see \Cref{specialcase}. A plausibility argument for an object-wise equivalence upon restricting to $\Glo_\ab$ is simple: given a fixed $n\geq 1$ and a geometric preoriented abelian sheaf $\FF$ over some affine scheme $\s=\Spec \A$, then the values of both
\[\Ell_{\FF/\s}^{\B C_n} \qquad \text{and} \qquad \Temp_{\FF[\P^\infty]/\s}^{\B C_n}\simeq \Spec \A_{\FF[\P^\infty]}^{\B C_n}\]
are both given by the $n$-torsion $\FF[n]$ of $\FF$; this fact is built into the constructions of both Gepner--Meier and Lurie.\\

Constructing the comparison maps of \Cref{maintheoremalg,maintheoremgeo} is not a trivial task. First, there is the formal issue of lifting Lurie and Gepner--Meier's constructions to functors out of a Grothendieck construction (\Cref{functorialell,functorialtemp}); both of these statements are hinted at in their original sources. The next problem is to choose what $\infty$-category everything is to take place in. Indeed, Lurie works with a functor-of-points perspective to study $\P$-divisible groups over \emph{affine} derived stacks, and his $\infty$-categories of orbispaces are based on finite abelian groups. On the other hand, Gepner--Meier work with a more geometric point-of-view to study elliptic curves over general derived stacks, and their $\infty$-category of orbispaces is based on arbitrary compact Lie groups. As neither approach is a strict generalisation of the other, we both extend Lurie's theory to general derived stacks $\s$ and geometrise it to fit into Gepner--Meier's language, as well as restrict Gepner--Meier's work to global sections and orbispaces on finite abelian groups to fit into Lurie's framework---hence the two flavours of comparison theorems above.

\subsection*{Application to equivariant topological modular forms}

As an application, we use \Cref{maintheoremgeo,maintheoremalg} to transport results between these two \emph{a priori} different theories. In particular, we apply them to the aforementioned $G$-fixed points of equivariant topological modular forms $\TMF^{\B G}$ for certain compact Lie groups $G$. The following is obtained by combining arguments of Gepner--Meier \cite[Rmk.10.6]{davidandlennart} and Lurie \cite{ec3} together with \Cref{maintheoremalg}; see \Cref{dualisableforfinitegroups}.

\begin{theoremalph}\label{dualisableproperty}
Let $G$ be a compact Lie group which can be written as the product of a torus with a finite group and $\x$ be an \textbf{oriented} elliptic curve over a derived stack $\s$. Then the $\E_\infty$-object $\O_{\Ell_{\x/\s}^{\B G}}$ in the $\infty$-category $\QCoh(\s)$ is dualisable. In particular, $\TMF^{\B G}$ is \emph{dualisable} as a $\TMF$-module.
\end{theoremalph}

This confirms a conjecture of Gepner--Meier for the above class of compact Lie groups; see \cite[Cor.1.4 \& Rmk.10.6]{davidandlennart}. There are other results from \cite{ec3} about tempered cohomology that one could try to pass to equivariant elliptic cohomology, such as an Atiyah--Segal completion theorem or a base-change result for $\pi$-finite spaces. In the affine case, \Cref{maintheoremalg} transfers these results seamlessly to equivariant elliptic cohomology. For nonaffine $\s$, one wants to commute a limit past various colimits, so the result does not immediately follow; we are exploring other methods to extend these results to the nonaffine setting in forthcoming work of the author and Balderrama--Linskens.

\subsection*{Conjectural refinement to global homotopy theory}

We believe that the comparison between tempered and equivariant elliptic cohomology can be further enhanced. Let us write $\Sp_\gl$ for the symmetric monoidal $\infty$-category of \emph{global spectra} of Schwede \cite{s}. By the main theorem of \cite{denissilluca}, a global spectrum $X$ can be viewed as a collection of $G$-equivariant spectra $X_G$ for each compact Lie group $G$, a collection of maps $f_\varphi \colon \varphi^\ast X_K \to X_G$ for each homomorphism of compact Lie groups $\varphi\colon G \to K$ which are equivalences if $\varphi$ is injective, together with a lot of higher coherences. There is also a version of global spectra based only on finite groups, which we denote as $\Sp_{\fin-\gl}$, and receives a restriction functor $\res_\fin \colon \Sp_\gl \to \Sp_{\fin-\gl}$. In \cite{gepner2024global2ringsgenuinerefinements}, Gepner--Linskens--Pol prove that given an \emph{oriented} geometric abelian sheaf, one can refine the global section functors $\Ga \Ell_{\FF/\s}, \Ga\Temp_{\FF[\P^\infty]/\s}\colon \OS^\op \to \CAlg_{\Ga(\s)}$ to $\E_\infty$-objects in $\Sp_\gl$ and $\Sp_{\fin-\gl}$, respectively. In work to appear together with Balderrama--Linskens \cite{tempered}, we provide an alternative to this Gepner--Linskens--Pol construction. Let us denote these global refinements by $\Ga_\gl \Ell_{\FF/\s}$ and $\Ga_{\fin-\gl}\Temp_{\FF[\P^\infty]/\s}$, respectively.\\

By construction, for a finite group $G$, the $G$-fixed points of these global refinements are the $\E_\infty$-rings $\Ga \Ell_{\FF/\s}^{\B G}$ and $\Ga\Temp_{\FF[\P^\infty]/\s}^{\B G}$, respectively. By \Cref{maintheoremgeo,maintheoremalg}, these $\E_\infty$-rings naturally agree. This leads us to the following:

\begin{conjecturealph}\label{globalhomotopytypesagree}
Let $\FF$ be a geometric oriented abelian sheaf over a derived stack $\s$. Then the two $\E_\infty$-objects in $\Sp_{\fin-\gl}$
\[\res_\fin \Ga_\gl \Ell_{\FF/\s} \qquad \text{and} \qquad \Ga_{\fin-\gl}\Temp_{\FF[\P^\infty]/\s}\]
are naturally equivalent.
\end{conjecturealph}

We believe this conjecture is true: as soon as one has a map of $\E_\infty$-objects between this pair of objects in $\Sp_{\fin-\gl}$, then it is most likely to be an equivalence as by \Cref{maintheoremgeo,maintheoremalg} these objects agree on fixed points. However, it does not follow from \cite{gepner2024global2ringsgenuinerefinements} that such a map exists. Indeed, there are two main difficulties. The first is that although it is simple to use \Cref{maintheoremgeo,maintheoremalg} to construct a map of naïve global $2$-rings (\cite[Df.6.1]{gepner2024global2ringsgenuinerefinements}) between those coming from equivariant elliptic cohomology and tempered cohomology, the associated genuine global $2$-rings are ``$\mathcal{T}$-genuine'' with respect to different families of groups $\mathcal{T}$; see \cite[Df.11.10]{gepner2024global2ringsgenuinerefinements}. Equivariant elliptic cohomology defines a tori-genuine global $2$-ring and tempered cohomology defines a finite-genuine global $2$-ring; see \cite[Ths.14.14 \& 15.24]{gepner2024global2ringsgenuinerefinements}. In fact, the former is explicitly \textbf{not} genuine with respect to finite groups, see \cite[Rmk.14.16]{gepner2024global2ringsgenuinerefinements}, and the latter is not even defined with respect to compact Lie groups. The way around this seems to be to develop a theory of \emph{tempered sheaves} which specialises to both quasi-coherent sheaves on tori and tempered local systems on finite groups; this is current work in progress by {Konovalov}--{Perunov}--{Prikhodko} \cite{temperedsheaves}.\\

Another question with this conjecture concerns how functorial the global refinements of \cite{gepner2024global2ringsgenuinerefinements} actually are. One can show that these global refinements are \textbf{not} functorial in the $\infty$-categories $\int_{\Stk} \PregeoAb$ appearing in \Cref{maintheoremgeo,maintheoremalg}. At its core, this is due to the fact that not all morphisms of preoriented abelian sheaves commute with \emph{transfers} on the associated equivariant cohomology theories. For example, Hirata--Kono show that Adams operations $\psi^k$ commute with transfers on equivariant topological $K$-theory if and only if $k$ is coprime to the group order; see \cite[Th.3.1]{hiratakono}. The functoriality of $\Ga_{\fin-\gl} \Temp_{\FF[\P^\infty]/\s}$ will be more thoroughly described in forthcoming work with Balderrama--Linskens \cite{tempered}.\\

We also expect this comparison to lift to higher equivariant multiplicative structures such as \emph{norms} and \emph{ultra-commutativity} \cite[\textsection 5]{s}, when such structures are conjured up from algebro-geometric input. This is further explored in \cite{normsontmf}.

%\addcontentsline{toc}{subsection}{Outline}
\subsection*{Outline}

In \Cref{geosheavessection}, we review Lurie's notions of abelian group objects and torsion objects. In \Cref{geometricsection}, we define \emph{geometric} abelian sheaves over a stack $\s$ after a little set-up, and prove some of their basic properties. In \Cref{homomorphismsheafsection}, we compare homomorphism sheaves into geometric abelian sheaves and their associated torsion objects. In \Cref{ecsection}, we recall and reiterate some basic facts about equivariant elliptic cohomology, drawing our references from \cite{davidandlennart}. In \Cref{tempcohomsection}, we similarly review Lurie's tempered cohomology, following \cite{ec3}, and discuss one possible extension from affine derived stacks to general derived stacks. We emphasise the functoriality inherent in these constructions of Gepner--Meier and Lurie. In \Cref{comparisonsection}, we connect our two different notions of preorientations and then piece everything together to prove \Cref{maintheoremalg,maintheoremgeo}. In \Cref{dualisableiltysection}, we use \Cref{maintheoremalg} to prove \Cref{dualisableproperty} by reducing to an argument of Gepner--Meier based on Lurie's tempered ambidexterity.

%\addcontentsline{toc}{subsection}{Notation}
\subsection*{Notation}

The language of $\infty$-categories will be used extensively throughout this article. In particular, we will call an $\infty$-category simply a \emph{category} and write $\Cat$ for the ($\infty$-) category of ($\infty$-) categories. The Cartesian unstraightening of a functor $F\colon \C^\op \to \Cat$ will be written as $\int_\C F \to \C$ and the coCartesian unstraightening of a functor $G\colon \D \to \Cat$ as $\int^\D G \to \D$. We will also use the following notation associated with categories of orbispaces:
\begin{itemize}
	\item $\Glo$ is the category of spaces of the form $\B G$ for some finite group $G$. We write $\Glo_\ab$ for the full subcategory of $\Glo$ spanned by those spaces with abelian $\pi_1$.
	\item $\Glo_{\Lie}$ is the global indexing category of compact Lie groups as defined in \cite[Df.2.7]{davidandlennart}, whose objects are the classifying stacks $\B G$ for compact Lie groups $G$ and morphisms spaces $\Map_{\Glo_\Lie}(\B G, \B K)$ given by $\Hom_{\Lie}(G,K)_{h K}$. This also has a full subcategory $\Glo_{\ab\Lie}$ spanned by those $\B G$ where $G$ is abelian; we will not use these categories other than for exposition.
	\item $\OS$ is the category of \emph{orbispaces}, so the presheaf category $\Fun(\Glo^\op,\Spc)$. We write $\OS_\ab$ for the category of \emph{abelian orbispaces}, so the presheaf category on $\Glo_\ab$, or equivalently, the full subcategory of $\OS$ generated under small colimits by representables of the form $\B G$ for a finite abelian group $G$.
	\item $\OS_\Lie$ is the category of \emph{compact Lie orbispaces}, so the presheaf category of $\Glo_{\Lie}$. This also has a subcategory $\OS_{\ab\Lie}$ defined in the obvious way; we will not use these categories other than for exposition.
\end{itemize}

In the main text of this paper, a \emph{stack} will mean a \emph{qcqs Noetherian nonconnective spectral Deligne--Mumford stack}, see \cite[\textsection1.4]{sag}, and a morphism of stacks will be almost of finite presentation (afp); see \cite[Df.5.1]{davidandlennart}. These restrictions are merely technical requirements to match with the big \'{e}tale site of Gepner--Meier and play no other serious role. In particular, the category of stacks $\Stk$ has only afp morphisms. For a stack $\s$, we will write $\Shv(\s)$ for the category of étale sheaves on the \emph{big \'{e}tale site} of $\s$, following \cite[Df.5.3 \& Lm.D.4]{davidandlennart}. In particular, the big étale site is defined as the category of stacks of almost finite presentation over $\tau_{\geq 0}\s$. If $\s=\Spec \A$ is affine, we will write $\Shv(\s)=\Shv(\A)$. If $\FF$ is a sheaf (or [pre]abelian group object) over $\s$ and $f\colon \s'\to \s$ is a morphism of stacks, we will write $f^\ast \FF$ for the base-change of $\FF$ to an object over $\s'$. We also abbreviate $\CAlg(\QCoh(\s))$ to $\CAlg(\s)$.

%\addcontentsline{toc}{subsection}{Acknowledgements}
\subsection*{Acknowledgements}
Thank you to William Balderrama, Sil Linskens, Lennart Meier, Lucas Piessevaux, and Julie Rasmusen for our discussions involving (non)equivariant elliptic and tempered cohomology as well as the surrounding mathematics. Thank you as well to Lennart and Sil for reading a draft. Finally, thank you to anonymous referee(s) for their numerous clarifying comments and helpful suggestions, all of which helped to improve this article. The author is, at the time of writing, an associate member of the Hausdorff Center for Mathematics at the University of Bonn (\texttt{DFG GZ 2047/1}, project ID \texttt{390685813}).

%%%%%%%%%%%%%%%%%%%%%%%%%%%%%%%%%%%%%%%%%%%%%%%%%%%%%%%%%%%
%%%%%%%%%%%%%%%%%%%%%%%%%%%%%%%%%%%%%%%%%%%%%%%%%%%%%%%%%%%
%%%%%%%%%%%%%%%%%%%%%%%%%%%%%%%%%%%%%%%%%%%%%%%%%%%%%%%%%%%

\section{Naturality of torsion objects}\label{geosheavessection}
Before we review and compare the works of Gepner--Meier and Lurie, we need to establish some common notation and conventions. In \cite{davidandlennart}, Gepner--Meier build equivariant elliptic cohomology out of abelian sheaves, and in \cite{ec3}, Lurie builds tempered cohomology out of $\P$-divisible groups. Our first goal is to study a certain subcategory of abelian sheaves whose associated \emph{torsion objects} are $\P$-divisible groups. We begin with a recollection of abelian group objects and torsion objects in general.

\begin{mydef}[{\cite[Df.1.2.4]{ec1}}]
	Let $\C$ be a category with finite limits. Write $\Ab_{\fgf}$ for the subcategory of the $1$-category of abelian groups spanned by those finitely generated free abelian groups. An \emph{abelian group object} in $\C$ is a product preserving functor from $\Ab_{\fgf}^\op$ to $\C$. We write $\Ab(\C)$ for the category of abelian group objects of $\C$ as a full subcategory of $\Fun(\Ab_{\fgf}^\op, \C)$. If $\C = \Shv(\s)$ for a stack $\s$, then we write $\Ab(\Shv(\s)) = \Ab(\s)$.
\end{mydef}

\begin{mydef}[{\cite[Df.6.4.2]{ec1}}]\label{torsionobjectsdef}
	Let $\C$ be a category with finite limits. Write $\Ab_\fin$ for the subcategory of the $1$-category of abelian groups spanned by finite abelian groups. A \emph{torsion object} of $\C$ is a functor $F\colon \Ab_\fin^\op \to \C$ which commutes with finite products and for every short exact sequence of finite abelian groups
	\[ 0 \to H' \to H \to H'' \to 0 ,\]
	the diagram
	\[\begin{tikzcd}
	{F(H'')}\ar[r]\ar[d]	&	{F(H)}\ar[d]	\\
	{F(0)}\ar[r]			&	{F(H')}
	\end{tikzcd}\]	
	is a pullback diagram in $\C$. We write $\Tors(\C)$ for the category of torsion objects in $\C$. If $\C=\Shv(\s)$ for a stack $\s$, we will simply write $\Tors(\s)$.
\end{mydef}

\begin{example}\label{exmapleoftorsion}
	There is a chosen equivalence of $\infty$-categories $\Ab(\Spc) \simeq \Mod_\Z^\cn$ of abelian group objects in spaces and connective $\Z$-modules; see \cite[Ex.1.2.9]{ec1}. Similarly, one can identify $\Tors(\Spc)$ as the full subcategory $\Mod_\Z^{\cn,\Tors}$ spanned by those connective $\Z$-modules $M$ such that $M\otimes_\Z \Q=0$; see \cite[Ex.6.4.11]{ec1}.
\end{example}

The following is a slight refinement of \cite[Pr.6.4.6]{ec1}: it states that one can naturally associate a torsion object to an abelian group object.

\begin{prop}\label{naturalityoftorsionfunctor}
There is a functor $[\P^\infty]\colon \Cat^{\lex}\to \Cat^{\Delta^1}$, where $\Cat^\lex$ is the category of categories admitting finite limits and functors which preserve such limits, whose value on a category $\C$ is the functor $[\P^\infty]_\C\colon \Ab(\C)\to \Tors(\C)$ of \cite[Pr.6.4.6]{ec1}. In other words, $[\P^\infty]$ is natural in finite limit preserving functors. Restricting to categories with finite limits and sequential colimits that are exact, then the functors $[\P^\infty]_\C$ all have fully faithful left adjoints.
\end{prop}

The subscript of $[\P^\infty]_\C$ will be dropped when the context is clear. We will often use the fully faithful left adjoint to $[\P^\infty]_\C$ to view $\Tors(\C)$ as a subcategory of $\Ab(\C)$ and hence view $[\P^\infty]_\C$ as an endofunctor.

\begin{proof}
The assignments sending a category $\C$ with finite limits to $\Ab(\C)$ or $\Tors(\C)$ both define functors
\[\Ab,\Tors\colon \Cat^{\lex}\to \Cat,\]
respectively. The functor $[\P^\infty]$ above will then be a natural transformation between these functors. Recall that $\Ab(\C)$ is a full subcategory of $\Fun(\Ab_{\fgf}^\op,\C)$ and $\Tors(\C)$ is a full subcategory of $\Fun(\Ab_{\fin}^{\op},\C)$; see \cite[\textsection6.4]{ec3}. Let us write $\Ab_{\fgf,\fin}$ for the full subcategory of the $1$-category of abelian groups spanned by groups that are either finite abelian or finitely generated free. Part of the content of \cite[Pr.6.4.6]{ec3} states that the composition
\[\Ab(\C)\subseteq \Fun(\Ab_{\fgf}^\op,\C) \xrightarrow{\Ran} \Fun(\Ab_{\fgf,\fin}^\op,\C) \xrightarrow{\res} \Fun(\Ab_{\fin}^\op,\C)\]
factors through $\Tors(\C)$, which Lurie defines to be $[\P^\infty]$. Clearly right Kan extension and restriction along the inclusion $\Ab_{\fin}\subseteq \Ab_{\fgf,\fin}$ are naturally in $\C$, hence $[\P^\infty]$ is natural in $\C$. If we further assume that $\C$ has sequential colimits which are exact, then \cite[Pr.6.4.9]{ec1} exactly states that the functor $[\P^\infty]_\C$ has a left adjoint.
\end{proof}

\begin{mydef}\label{counitdefinition}
Given a category $\C$ with finite limits and sequential colimits that are exact, then we write $\varepsilon_\C\colon [\P^\infty]_\C\to \id_{\Ab(\C)}$ for the natural transformation of endofunctors of $\Ab(\C)$, the counit of the adjunction from \Cref{naturalityoftorsionfunctor}---we omit notation for the fully faithful left adjoint. The formation of this counit is natural in functors which preserve finite limits by \Cref{naturalityoftorsionfunctor}, meaning for such a functor $F\colon \C\to \D$, there is a natural equivalence $F(\epsilon_\C)\simeq \varepsilon _\D$ of natural transformation of endofunctors of $\Ab(\D)$. To be more precise, \Cref{naturalityoftorsionfunctor} allows us to define a functor
\[\Cat^{\lex+}\to \Cat \qquad \C\mapsto \Fun(\Ab(\C),\Ab(\C))_{[\P^\infty]_\C//\id_{\Ab(\C)}}\]
where $\Cat^{\lex+}$ is the full subcategory of $\Cat^\lex$ spanned by categories with all sequential colimits that are also exact. The naturality of the counit construction can then be described by defining a section of the associated coCartesian fibration:
\[\Cat^{\lex+}\to \int^{\Cat^{\lex+}} \Fun(\Ab(-),\Ab(-))_{[\P^\infty]_{(-)}//\id_{\Ab(-)}}\qquad \C\mapsto (\epsilon_\C\colon [\P^\infty]_\C\to \id_{\Ab(\C)}).\]
\end{mydef}

%%%%%%%%%%%%%%%%%%%%%%%%%%%%%%%%%%%%%%%%%%%%%%%%%%%%%%%%%%%
%%%%%%%%%%%%%%%%%%%%%%%%%%%%%%%%%%%%%%%%%%%%%%%%%%%%%%%%%%%
%%%%%%%%%%%%%%%%%%%%%%%%%%%%%%%%%%%%%%%%%%%%%%%%%%%%%%%%%%%

\section{Geometric abelian sheaves}\label{geometricsection}
Now we get into some spectral algebraic geometry. The goal of this section is to define \emph{geometric} abelian sheaves. These will be abelian sheaves over a stack whose associated torsion objects form a \emph{$\P$-divisible group}.\\

First, we need the following well-known lemmata; both of these facts are also recorded in more generality in \cite[Cor.3.4.4.4 \& Cor.2.1.2.4]{reconstruction}. Write $\Stk_{/\s}^\flat$ for the subcategory of $\Stk_{/\s}$ spanned by maps of stacks $\t \to \s$ which are flat.

\begin{lemma}\label{flatnessandconnectivecovers}
	Let $\s$ be a stack. Then the connective cover functor
	\[\tau_{\geq 0} \colon\Stk_{/\s}^{\flat} \to \Stk_{/\tau_{\geq 0}\s}^{\flat}\]
	is an equivalence.
\end{lemma}

\begin{proof}
	The inverse to the connective cover functor is given by base-change:
\[\Stk_{/\tau_{\geq 0}\s}^{\flat}\to \Stk_{/\s}^{\flat}\qquad (\x_0\to \tau_{\geq 0}\s)\mapsto (\x_0\times_{\tau_{\geq 0}\s} \s\to \s).\]
	Checking that these functors are inverse to each other can be done étale locally, so it suffices to check this on affines. This then follows from \cite[Pr.7.2.2.24]{ha}; essentially a Tor-spectral sequence concentrated in filtration $0$.
\end{proof}

\begin{lemma}\label{equivalencewithconnectiveguys}
Let $\s$ be a stack. Then the functor
\[\Stk_{/\s}\to \Shv(\s)\qquad (\x\to \s)\mapsto \left((\t\to \tau_{\geq 0}\s)\mapsto \Map_{\Stk_{/\tau_{\geq 0}\s}}(\t,\tau_{\geq 0}\x)\right)\]
is fully faithful when restricted to $\Stk_{/\s}^{\afp,\flat}$, so those morphisms of stacks $\x\to\s$ which are flat and almost of finite presentation.
\end{lemma}

This functor is not fully faithful with these restrictions. Without the flatness hypothesis, the connective cover functor $\Stk_{/\s} \to \Stk_{/\tau_{\geq 0}\s}$ is far from an equivalence or even being fully faithful. Without the almost of finite presentation hypothesis, the functor $\Stk_{/\tau_{\geq 0}\s} \to \Shv(\s)$ sending a stack to the functor it represents is also not clearly fully faithful.

\begin{proof}
The given functor, hence also its restriction, can be factored as the composite
\[\Stk_{/\s}\xrightarrow{\tau_{\geq 0}} \Stk_{/\tau_{\geq 0}\s}\to \Shv(\s)\]
where the second map is now given by a variant of the Yoneda embedding---its restriction to $\Stk^{\afp}_{/\tau_{\geq 0}\s}$ is the classical Yoneda embedding. It now suffices to see that the induced functor
\[\Stk_{/\s}^{\afp,\flat}\xrightarrow{\tau_{\geq 0}} \Stk_{/\tau_{\geq 0}\s}^{\afp,\flat}\]
is fully faithful, which follows from \Cref{flatnessandconnectivecovers}.
\end{proof}

Given a stack $\s$ and an abelian $\s$-sheaf $\FF$, we can use the functor $[\P^\infty]_{\Shv(\s)}$ to obtain a torsion object $\FF[\P^\infty]\colon \Ab_\fin^\op \to \Shv(\s)$. We now need to define what it means for this torsion object to be a \emph{$\P$-divisible group} over $\s$. Let us start with the affine case; see \cite[Df.2.6.1]{ec3} for this definition.

\begin{mydef}\label{pdivdef}
A $\P$-divisible group $\G$ over an $\E_\infty$-ring $\A$ is a functor $\G\colon \CAlg^\cn_{\tau_{\geq 0}\A} \to \Mod_\Z^\cn$ such that for each $\mathsf{B}$, the $\Z$-module $\G(\mathsf{B})$ is torsion, so $\G(\mathsf{B}) \otimes_\Z \Q = 0$, for every finite abelian group $H$ the functor
\[(\mathsf{B}\in \CAlg_{\tau_{\geq 0\A}}^\cn) \mapsto (\Map_{\Mod_\Z^\cn}(H,\G(\mathsf{B}))\in \Spc)\]
is corepresented by a finite flat $\E_\infty$-$\tau_{\geq 0} \A$-algebra $\O_{\G[H]}$, and for every positive integer $n$, the map $[n]\colon \G \to \G$ is locally surjective in the finite flat topology. Write $\BT(\A)$ for the category of $\P$-divisible groups over $\A$, defined as the full subcategory of $\Fun(\CAlg^\cn_{\tau_{\geq 0}\A}, \Mod_\Z^\cn)$ spanned by $\P$-divisible groups.
\end{mydef}

This construction of categories of $\P$-divisible groups over $\E_\infty$-rings is functorial in $\E_\infty$-rings.

\begin{mydef}\label{categoriesofpdivisiblegroups}
Given a morphism $f\colon \A \to \A'$ and a $\P$-divisible $\G \colon \CAlg_\A \to \Mod_\Z^\cn$, we write $f^\ast \G$ for the $\P$-divisible group given by the composite of $\G$ with the forgetful functor
\[\CAlg_{\A'} \xrightarrow{f_\ast} \CAlg_\A \xrightarrow{\G} \Mod_\Z^\cn.\]
By \cite[Rmk.6.5.4]{ec1}, this defines a functor $\BT\colon \CAlg\to\Cat$ sending $\A$ to $\BT(\A)$.
\end{mydef}

\begin{remark}\label{remarkalternativedefinition}
	Let $\A$ be an $\E_\infty$-ring. By \cite[Pr.6.5.5]{ec1}, the functor
	\[\BT(\A) \to \Fun((\Ab_\fin^\op, \Stk_{/\Spec \A}), \qquad \G \mapsto (H \mapsto \O_{\G[H]}),\]
	where $\O_{\G[H]}$ is the finite flat $\E_\infty$-$\A$-algebra from \Cref{pdivdef}, is fully faithful, with essential image those functors $F$ which are torsion objects in $\Stk_{/\Spec \A}$ and such that for each prime $p$ and each positive integer $n\geq 0$, the map $F(\Z/p^n\Z) \to F(p\Z/p^n\Z)$ is a finite flat surjection.
\end{remark}

Our generalisation to $\P$-divisible groups over stacks is totally formal.

\begin{mydef}
Let $\M_{\BT}$ be the \emph{moduli functor of $\P$-divisible groups}, so the functor
\[\M_{\BT}\colon \CAlg^\cn\to \Spc\qquad \A\mapsto \BT(\A)^\simeq\]
given by taking the groupoid core of $\BT(\A)$ from \Cref{categoriesofpdivisiblegroups}. A \emph{$\P$-divisible group over a stack $\s$} is a natural transformation of functors
\[\G\colon \tau_{\geq 0}\s\Longrightarrow \M_\BT\colon \CAlg^\cn\to \Spc,\]
so $\M_{\BT}(\s) = \Map(\tau_{\geq 0}\s,\M_\BT)$. Equivalently, writing $\s$ as a colimit of representables $\s = \colim\Spec \A$ in the big étale site of $\s$, we see the natural map
\[\M_{\BT}(\s)\xrightarrow{\simeq} \lim_{\s_\et^\aff} \M_{\BT}(\A)\]
into the limit taking over all étale morphisms $f\colon \Spec \A \to \s$ from an affine,\footnote{The category $\s_\et^\aff$ has many useful cofinal subdiagrams. Indeed, we assumed that our stacks $\s$ are qcqs, so by \cite[A.17]{luriestheorem}, $\s$ admits an affine \'{e}tale hypercover $\Spec \A_\bullet\to \s$ which can be used to calculate limits over $\s_\et^\aff$.} is an equivalence; see \cite[Pr.3.2.2(5)]{ec2}. In fact, the proof of \cite[Pr.3.2.2(5)]{ec2} suggests that defining the category of $\P$-divisible groups over $\s$ as the limit
\[\BT(\s) = \lim_{\s_\et^\aff} \BT(\A)\]
is reasonable. Writing $\Mod_\Z^{\cn,\Tors}$ for the subcategory of $\Mod_\Z^\cn$ spanned by torsion $\Z$-modules, as in \Cref{exmapleoftorsion}, define the value of a $\P$-divisible group $\G$ over a stack $\s$ as the limit inside $\Mod_\Z^{\cn,\Tors}$
\begin{equation}\label{globalsectionsofpdiv}\G(\s)=\lim f^\ast\G(\A).\end{equation}
As maps of $\P$-divisible groups $\G \to \G'$ over $\s$ induce, by definition, compatible maps $f^\ast \G\to f^\ast \G'$ for each étale map $f\colon \Spec \A \to \s$, we obtain an induced map of $\Z$-modules $\G(\s) \to \G'(\s)$. In particular, this induces a functor $\BT(\s) \to \Mod_\Z^{\cn, \Tors}$.
\end{mydef}

\begin{warn}\label{pdivnosheaffootnote}
A $\P$-divisible group defines an fpqc-sheaf when viewed to take values in $\Mod_\Z^{\cn,\Tors}$, but not necessarily when viewed as taking values in $\Mod_\Z^\cn$; see \cite[Warn.6.5.10]{ec1}.
\end{warn}

	There are other slight variations on generalising Lurie's tempered cohomology to stacks, using a generalisation of \Cref{remarkalternativedefinition}, for example; this is our approach in forthcoming work together with Balderrama--Linskens on viewing tempered cohomology theories as generalised global homotopy theories. We prefer the above definition in this article though, as it simplifies our discussion of preorientations, which will appear in \Cref{defofpreorientation}.\\
	
	We can finally define the main character of this section: \emph{geometric abelian sheaves}.

\begin{mydef}\label{torsionofsheavesdefintiion}
	Let $\s$ be a stack and $\FF$ an abelian $\s$-sheaf. We say that $\FF$ is \emph{geometric} if
	\begin{enumerate} 
	\item	the associated torsion object $\FF[\P^\infty]\colon \Ab_\fin^\op \to \Shv(\s)$ factors through $\Stk^{\flat,\afp}_{/\s}$, using \Cref{equivalencewithconnectiveguys}, hence defines a torsion object of $\Stk_{/\s}$, and
	\item that for all étale maps $f\colon \Spec \A \to \s$, the torsion object $(f^\ast \FF)[\P^\infty] \colon \Ab_\fin^\op \to \Stk_{/\Spec \A}$ defines a $\P$-divisible group over $\A$ à la \Cref{remarkalternativedefinition}.
	\end{enumerate}
	We write $\geoAb(\s)$ for the full subcategory of $\Ab(\s)$ spanned by geometric abelian sheaves $\FF$.
\end{mydef}

\begin{remark}\label{geoabgivespdiv}
In the situation of \Cref{torsionofsheavesdefintiion}, an application of \Cref{naturalityoftorsionfunctor} to the exact functor $f^\ast \colon \Stk_{/\s} \to \Stk_{/\Spec \A}$ given by base-change $-\times_\s \Spec \A$ gives a natural equivalence $(f^\ast \FF)[\P^\infty] \simeq f^\ast(\FF[\P^\infty])$. In particular, the collection of $\P$-divisible groups $f^\ast(\FF[\P^\infty])$ associated with a geometric abelian $\s$-sheaf $\FF$ define a $\P$-divisible group $f^\ast (\FF[\P^\infty])$ over $\s$.
\end{remark}

\begin{example}
	By \cite[Pr.6.7.1]{ec1}, the sheaf represented by an abelian variety $\x$ over a stack $\s$ is a geometric abelian $\s$-sheaf. This encompasses elliptic curves over a stack $\s$ by \cite[Df.2.0.2]{ec1}. Similarly, abelian $\s$-stacks which are afp and flat, as used in \cite[Con.14.12]{gepner2024global2ringsgenuinerefinements} are also abelian $\s$-sheaves by \Cref{equivalencewithconnectiveguys}, which are further geometric abelian $\s$-sheaves if the underlying torsion object in $\s$-stacks defines a $\P$-divisible group over $\s$.
\end{example}

The assignment sending a geometric abelian $\s$-sheaf to a $\P$-divisible group is functorial:

\begin{prop}\label{firstattemptatfunctorusingpinfty}
	The assignment sending a geometric abelian $\s$-sheaf $\FF$ to the $\P$-divisible group $\FF[\P^\infty]$ over $\s$ defines a morphism of Cartesian fibrations over $\Stk$
\[[\P^\infty] \colon \int_{\Stk} \geoAb \to \int_{\Stk}\BT.\]
\end{prop}

\begin{proof}
	By straightening/unstraightening, it suffices to construct a natural transformation of functors from $\Stk^\op \to \Cat$ from $\geoAb(-)$ to $\BT(-)$. On affines $\Spec \A$, this natural transformation is given by first considering $[\P^\infty]_{\Shv(\Spec \A)} \colon \Ab(\Spec \A) \to \Tors(\Spec \A)$ of \Cref{naturalityoftorsionfunctor}. Restricting this functor to $\geoAb(\Spec \A)$ then factors through $\Tors(\Stk^{\flat,\afp}_{/\Spec \A})$ by definition, and further through $\BT(\A)$ by \Cref{remarkalternativedefinition}. This is the desired natural transformation on affines. For a general stack $\s$, one maps into the limit $\BT(\s) = \lim \BT(\A)$ over the small étale site of $\s$.
\end{proof}

The notation $\FF[\P^\infty]$ used above is potentially ambiguous.

\begin{warn}\label{potentialpdivnotbeingsheaves}
	Given a geometric abelian $\s$-sheaf $\FF$, we will write $\FF[\P^\infty]$ for two different, but related, objects: the first, is as an object in $\Tors(\s)$ given by applying $[\P^\infty]_{\Shv(\s)} \colon \Ab(\s) \to \Tors(\s)$ to $\FF$ as in \Cref{naturalityoftorsionfunctor}. As $\Shv(\s)$ is complete and cocomplete and sequential colimits are exact, we can also apply the fully faithful left adjoint to $[\P^\infty]_{\Shv(\s)}$ and consider $\FF[\P^\infty]$ as an object in $\Ab(\s)$. In particular, the natural counit map of \Cref{counitdefinition} gives a comparison map $\FF[\P^\infty] \to \FF$. On the other hand, by \Cref{geoabgivespdiv}, we can also view $\FF[\P^\infty]$ as a $\P$-divisible group over $\s$, which is not \emph{a priori} an $\s$-sheaf. 
\end{warn}

The above warning will not cause us any problems as we are really interested in certain \emph{homomorphism stacks}.

\begin{mydef}\label{defhomomorphissheaf}
Let $\s$ be a stack. Write $\HOM_\s(\mathscr{F},\mathscr{G})$ for the internal mapping object $\MAP_{\Ab(\s)}(\mathscr{F},\mathscr{G})$ between two abelian $\s$-sheaves $\mathscr{F}$ and $\mathscr{G}$.\footnote{Recall that we are writing $\Ab(\s)$ as an abbreviation of $\Ab(\Shv(\s))$, so abelian group objects in sheaves on the big étale site of $\s$. This is equivalent to the category of sheaves on the big étale site of $\s$ taking values in abelian group objects in spaces. In particular, the homomorphism sheaf $\HOM_{\s}(\mathscr{F}, \mathscr{G})$ in $\Ab(\Shv(\s))$ is equivalent to the internal mapping sheaf in the category $\Shv_{\Ab(\Spc)}(\s)$.}
\end{mydef}

Let us write $\underline{(-)}\colon \Ab(\Spc) \to \Ab(\s)$ for the constant sheaf functor, the left adjoint to evaluation on $\s$.
\begin{mydef}\label{defhomsheafforpdiv}
	Let $\A$ be an $\E_\infty$-ring, $\G$ a $\P$-divisible group over $\A$, and $H$ a finite abelian group. We define $\HOM_\A(\underline{H}^\vee, \G)$ in $\Shv(\Spec \A)$ by the assignment
	\[(f\colon \t \to \Spec \tau_{\geq 0} \A) \mapsto \Map_{\Mod_\Z^\cn}(\underline{H}^\vee, f^\ast \G(\t)).\]
	This assignment is a sheaf as it is represented by the $\Spec \tau_{\geq 0}\A$-stack $\Spec \O_{\G[H^\vee]}$ from \Cref{pdivdef}---it is crucial here that $H^\vee$ is a finite abelian group, ie, lies in the subcategory $\Mod_\Z^{\cn,\Tors}$ of $\Mod_\Z^\cn$. This sheaf satisfies base-change along morphisms of $\E_\infty$-rings $\A \to \A'$ essentially from the definition of $\O_{\G[H^\vee]}$ as corepresenting the above functor and the universal property of base-change. For a general stack $\s$ and a $\P$-divisible group $\G$ over $\s$, we define the sheaf $\HOM_\s(\underline{H}^\vee, \G)$ as the object of 
	\[\Shv(\s) \simeq \lim_{\Spec \A \in \s_\et^\aff} \Shv(\Spec \A)\]
	given by the collection of the sheaves $\HOM_\A(\underline{H}^\vee, f^\ast\G)$ for $f\colon \Spec \A \to \s$.
\end{mydef}

Finally, we reconcile \Cref{potentialpdivnotbeingsheaves} in our cases of interest.

\begin{prop}\label{homoforpdivagreeswithtorsionongeometricguys}
	Let $\s$ be a stack, $\FF$ be a geometric abelian $\s$-sheaf, and $G$ a finite abelian group. Then the two definitions of $\HOM_\s(\underline{G}^\vee, \FF[\P^\infty])$, considering $\FF[\P^\infty]$ either as an abelian $\s$-sheaf (\Cref{defhomomorphissheaf}) or a $\P$-divisible group (\Cref{defhomsheafforpdiv}), naturally agree.
\end{prop}

\begin{proof}
	As-per-usual, this equivalence of abelian $\s$-sheaves for a general stack $\s$ will be a limit of such sheaves over affine stacks $\Spec \A$, so we are reduced to the affine case of $\s = \Spec \A$. Let $f\colon \t \to \Spec \tau_{\geq 0}\A$ be an $\Spec \tau_{\geq 0}\A$-stack. Then considering $\FF[\P^\infty]$ as an abelian $\Spec \A$-sheaf, we have
	\[\HOM_{\Spec \A}(\underline{G}^\vee, \FF[\P^\infty])(\t) = \Map_{\Ab(\t)}(\underline{G}^\vee, f^\ast(\FF[\P^\infty])) \simeq \Map_{\Ab(\Spc)}(G^\vee, f^\ast(\FF[\P^\infty])(\t)),\]
	using only the natural constant $\t$-sheaf-evaluation adjunction. The right-hand side is precisely the value of $\HOM_{\Spec \A}(\underline{G}^\vee, \FF[\P^\infty])$ on $\t$ considering $\FF[\P^\infty]$ as a $\P$-divisible group by definition; see \Cref{exmapleoftorsion} and \Cref{defhomsheafforpdiv}. As in \Cref{defhomsheafforpdiv}, it is crucial here that $G$ lies in the essential image of the embedding $\Tors(\Spc) \to \Ab(\Spc)$, which is clear as $G$ is finite.
\end{proof}

%%%%%%%%%%%%%%%%%%%%%%%%%%%%%%%%%%%%%%%%%%%%%%%%%%%%%%%%%%%
%%%%%%%%%%%%%%%%%%%%%%%%%%%%%%%%%%%%%%%%%%%%%%%%%%%%%%%%%%%
%%%%%%%%%%%%%%%%%%%%%%%%%%%%%%%%%%%%%%%%%%%%%%%%%%%%%%%%%%%
\section{Comparison of homomorphism sheaves}\label{homomorphismsheafsection}

For a stack $\s$, a geometric abelian $\s$-sheaf $\FF$, and a finite abelian group $G$, we now have two homomorphism sheaves
\[\HOM_\s(\underline{G}^\vee,\FF), \qquad \HOM_\s(\underline{G}^\vee,\FF[\P^\infty]).\]
We claim that the natural map $\FF[\P^\infty] \to \FF$ induces an equivalence between these two sheaves. We can boil this down to studying the associated \emph{$n$-torsion} sheaf for each $n$.

\begin{mydef}\label{def:ntorsion}
For any object $\GG\in \Ab(\s)$, let us write $\GG[n]$ for the \emph{$n$-torsion in $\GG$}, so the object of $\Ab(\s)$ defined by the Cartesian diagram
\[\begin{tikzcd}
{\GG[n]}\ar[r]\ar[d]	&	{\GG}\ar[d, "{[n]}"]	\\
{\s}\ar[r]		&	{\GG}
\end{tikzcd}\]
inside $\Ab(\s)$, where $[n]\colon \GG\to\GG$ is the $n$-fold multiplication map. Equivalently, $\GG[n]$ is the internal homomorphism sheaf $\HOM_{\s}(\underline{\Z/n\Z},\GG)$. A $\P$-divisible group $\G$ does not necessarily define a sheaf, see \Cref{pdivnosheaffootnote}, so in this case we \textbf{define} $\G[n]$ as the homomorphism sheaf $\HOM_\s(\underline{\Z/n\Z}, \G)$ of \Cref{defhomsheafforpdiv}, in other words, $\Spec \O_{\G[\Z/n\Z]}$.
\end{mydef} 

The special case of \Cref{homoforpdivagreeswithtorsionongeometricguys} when $G=\Z/n\Z$ then shows that the two natural definitions of the abelian $\s$-sheaf of $n$-torsion of $\FF[\P^\infty]$ agree. In fact, the $n$-torsion of $\FF[\P^\infty]$ also agrees with the $n$-torsion of $\FF$.

\begin{prop}\label{torsionsagreelemma}
Let $\FF$ be a geometric abelian $\s$-sheaf and fix a positive integer $n\geq 1$. Then the natural map $\FF[\P^\infty]\to \FF$ in $\Ab(\s)$, the counit of the adjunction from \Cref{counitdefinition}, induces an equivalence of finite flat abelian group stacks over $\s$
\[(\FF[\P^\infty])[n]\xrightarrow{\simeq} \FF[n].\]
\end{prop}

\begin{proof}
	From our definitions, it suffices to prove this in the affine case $\s = \Spec \A$. By \cite[Rmk.6.4.8]{ec1}, given an abelian group object $X$ in a category $\C$, the value of the torsion object $X[\P^\infty]$ on $\Z/n\Z$ is precisely the fibre of $[n] \colon X \to X$ in $\Ab(\C)$. For us, this means that we have the natural identification $\FF[\P^\infty](\Z/n\Z) = \FF[n]$. The value of $\FF[\P^\infty](\Z/n\Z)$ is given by the finite flat $\Spec \tau_{\geq 0}\A$-stack $\Spec \O_{\Z/n\Z}$ of \Cref{pdivdef} by definition. It follows from \Cref{defhomsheafforpdiv} that the desired sheaf of homomorphism from $\underline{\Z/n\Z}$ into the $\P$-divisible group $\FF[\P^\infty]$ is represented by this stack $\FF[\P^\infty](\Z/n\Z)$. By \Cref{homoforpdivagreeswithtorsionongeometricguys}, this is the same $(\FF[\P^\infty])[n]$.
\end{proof}

Note that the natural map $\FF[n]\to \FF$ factors through $\FF[\P^\infty]$, as $\FF[n]$ is a torsion object.

\begin{prop}\label{factoringhomstackslemma}
Let $\FF$ be a geometric abelian $\s$-sheaf over a stack $\s$ and $G$ a finite abelian group of order $|G|$. Then the natural maps
\[\FF[|G|]\to \FF[\P^\infty]\to \FF\qquad \in \Ab(\s)\]
induce equivalences of homomorphism sheaves over $\s$
\[\HOM_\s(\underline{G}^\vee,\FF[|G|])\xrightarrow{\simeq} \HOM_\s(\underline{G}^\vee,\FF[\P^\infty])\xrightarrow{\simeq} \HOM_\s(\underline{G}^\vee,\FF).\]
In particular, such homomorphism sheaves are represented by finite flat stacks over $\s$.
\end{prop}

\begin{proof}
Recall that the homomorphism stack $\HOM_\s(\underline{G}^\vee,\GG)$ in $\Shv(\s)$, for some $\GG\in \Ab(\s)$, is defined by sending $f\colon \t\to\tau_{\geq 0}\s$ of almost finite presentation to the space
\[\Map_{\Ab(\t)}(\underline{G}^\vee, f^\ast \GG).\]

An element in this mapping space, a homomorphism $\varphi\colon \underline{G}^\vee \to f^\ast \GG$, factors uniquely through the $|G|$-torsion of $f^\ast \GG$, as the $|G|$-fold multiplication map $[|G|]$ on $\underline{G}^\vee$ is zero:
\[\begin{tikzcd}
{\underline{G}^\vee}\ar[rrd, bend left = 30, "{\varphi}"]\ar[dr, dashed]\ar[rdd, bend right = 30]	&&\\
	&	{f^\ast \GG[|G|]}\ar[d]\ar[r]	&	{f^\ast \GG}\ar[d, "{[|G|]}"]	\\
	&	{\t}\ar[r]			&	{f^\ast \GG}
\end{tikzcd}\]

In other words, for any $\GG$ in $\Ab(\s)$, the natural map
\[\HOM_\s(\underline{G}^\vee,\GG[|G|]) \xrightarrow{\simeq} \HOM_\s(\underline{G}^\vee,\GG)\]
induced by the canonical inclusion $\GG[|G|]\to \GG$, is an equivalence inside $\Shv(\s)$. As the $|G|$-torsion of both $\FF[\P^\infty]$ and $\FF$ agree by \Cref{torsionsagreelemma}, we obtain the desired equivalence. To see that these homomorphism sheaves are representable by finite flat stacks over $\s$, we write $G^\vee\simeq \bigoplus \Z/n\Z$ as a sum of cyclic groups, which yields the equivalences of sheaves over $\s$
\[ \HOM_\s(\underline{G}^\vee,\FF) \simeq  \prod \HOM_\s(\underline{\Z/n\Z},\FF)\simeq \prod \FF[n].\]
Each $\FF[n]$ is finite flat over $\s$ by \Cref{torsionsagreelemma}, so we are done.
\end{proof}

This finishes our general theory. Now we can move on to the equivariant elliptic and tempered cohomology constructions of Gepner--Meier and Lurie.

%%%%%%%%%%%%%%%%%%%%%%%%%%%%%%%%%%%%%%%%%%%%%%%%%%%%%%%%%%%
%%%%%%%%%%%%%%%%%%%%%%%%%%%%%%%%%%%%%%%%%%%%%%%%%%%%%%%%%%%
%%%%%%%%%%%%%%%%%%%%%%%%%%%%%%%%%%%%%%%%%%%%%%%%%%%%%%%%%%%

\section{Functoriality of equivariant elliptic cohomology}\label{ecsection}

In \cite[Cons.3.13 \& 6.1]{davidandlennart}, Gepner--Meier construct a functor
\[ \Ell^{\GM}_{\FF/\s}\colon \OS_{\Lie}\to \Shv(\s)\]
for each preoriented abelian sheaf $\FF$ over a stack $\s$, which is functorial in the pair $\FF/\s$. We will remove the $\GM$-superscript shortly when we restrict to finite groups; see \Cref{definitionofellgeometric}. The first step of their construction is to consider \emph{preoriented abelian group objects} inside a category $\C$ with finite limits; see \cite[\textsection3]{davidandlennart}. 

\begin{notation}
	The representable objects of $\OS_\Lie$ are written as $\B G$ for some compact Lie group $G$. The space $B\T = B^2\Z$ can be viewed as an abelian group object in spaces either as a topological abelian group via \cite[Ex.3.3]{davidandlennart} or as $\Z[2]$ through the equivalence of categories $\Ab(\Spc) \simeq \Mod_\Z^\cn$ seen in \cite[Ex.1.2.9]{ec1}. Courtesy of \cite[Ex.3.4]{davidandlennart}, we can also recognise $\B G$ as an abelian group object of $\Glo_{\ab\Lie}$ (hence also $\OS_\Lie$) for each \emph{abelian} compact Lie group $G$.
\end{notation}
	
\begin{mydef}\label{generaldefintiionofpreorirnation}
	Let $\C$ be a category with finite limits and a terminal object $\ast$. A \emph{preoriented abelian group object in $\C$} is an object $\x$ in $\Ab(\C)$ equipped with a map $B \T\to \Map_\C(\ast,\x)$ of abelian group objects in spaces, where $\T=S^1$ is the circle. In other words, we define a category
\[\PreAb(\C)=\Ab(\C)\underset{\Ab(\Spc)}{\times} \Ab(\Spc)_{B \T/},\]
where $\Ab(\C) \to \Ab(\Spc)$ is corepresented by the terminal object $\ast$ in $\C$. When $\C = \Shv(\s)$ for a stack $\s$, we write $\PreAb(\s)$.
\end{mydef}

Another ingredient in their construction is the \emph{shifted Pontryagin dual}, defined as the functor
\[\Map_{\Glo_{\Lie}}(-,\B \T)\colon \Glo_{\ab\Lie}^\op\to \Ab(\Spc)_{B \T/}=\PreAb(\Spc),\qquad \B G\mapsto \widehat{\B G} = \Map_{\Glo_{\Lie}}(\B G,\B \T);\]
the abelian group structure on $\widehat{\B G}$ comes from that on $\B\T$ and the canonical preorientation on $\widehat{\B G}$ is induced by $\B G \to \ast$ and the equivalence $\widehat{\ast} \simeq B \T$. As shown in \cite[Con.3.8]{davidandlennart}, there is a natural identification
\begin{equation}\label{decompositionofdualofbg}
\widehat{\B G}\simeq B \T\times G^\vee
\end{equation}
in $\Ab(\Spc)$, where $G^\vee$ now denotes the classical Pontryagin dual of $G$, so the space of group homomorphisms into $\T$. One then constructs a functor
\begin{equation}\label{functorgivingellipticcohom}\PreAb(\s)\to \Fun(\Glo_{\ab\Lie},\Shv(\s))\qquad \y\mapsto \left(\B G\mapsto \MAP_{\PreAb(\s)}(\underline{\widehat{\B G}},\y)\right)\end{equation}
where $\MAP$ denotes the internal mapping object in $\PreAb(\s)$. Notice that the constant sheaf functor refines to a functor $\underline{(-)}\colon \PreAb(\Spc) \to \PreAb(\s)$, as for each object $B\T \to X$ of $\PreAb(\Spc)$, the unit $B\T \to \underline{B\T}(\s)$ provides us with a canonical preorientation of $\underline{X}$; see \cite[Con.3.13]{davidandlennart}.\\

For a preoriented abelian sheaf $\FF$ over a stack $\s$, the functor $\Ell^{\GM}_{\FF/\s} \colon \OS_\Lie \to \Shv(\s)$ is given by first left Kan extending the image of $\FF$ under the functor (\ref{functorgivingellipticcohom}) from $\Glo_{\ab\Lie}$ to $\Glo_{\Lie}$, and then a further left Kan extension along the Yoneda embedding from $\Glo_{\Lie}$ to $\OS_{\Lie}$. To encode the functoriality of this construction, we have the following statement. As $\Ab(-)$ is functorial in functors preserving finite limits, then by construction so is $\PreAb(-)$. In particular, each $f\colon \s'\to \s$ in $\Stk$ induces an exact pullback functor $f^\ast \colon \Shv(\s)\to \Shv(\s')$. The functoriality of $\Shv(-)$ and $\PreAb(-)$ then yield two functors
\[\PreAb\colon \Stk^\op\to \Cat,\qquad \s\mapsto \PreAb(\s)\]
\[\Fun^L(\OS_\Lie,\Shv(-))\colon \Stk^\op\to \Cat,\qquad \s\mapsto \Fun^L(\OS_\Lie,\Shv(\s)),\]
where $\Fun^L(-,-)$ denotes colimit preserving functors, whose associated Cartesian fibrations we write as $\int_{\Stk}\PreAb$ and $\int_{\Stk}\Fun^L(\OS_\Lie,\Shv(-))$.

\begin{prop}\label{functorialell}
There is a functor
\[\Ell^\GM\colon \int_{\Stk}\PreAb\to \int_{\Stk}\Fun^L(\OS_\Lie,\Shv(-))\]
whose value on a particular pair $\FF/\s$ in $\int_{\Stk}\PreAb$ is $\Ell^\GM_{\FF/\s}$ of Gepner--Meier.
\end{prop}

\begin{proof}
We will define our desired functor as the composite
\[
\Ell^\GM\colon \int_{\Stk} \PreAb \xrightarrow{\Ell^\GM|_{\Glo_{\ab\Lie}}} \int_{\Stk} \Fun(\Glo_{\ab\Lie}, \Shv(-)) \xrightarrow{\mathrm{Lan}}	\int_{\Stk} \Fun^L(\OS_\Lie,\Shv(-))
\]
where the second functor takes left Kan extensions along the inclusions $\Glo_{\ab\Lie}\to \Glo_\Lie$ and $\Glo_\Lie \to \OS_\Lie$. By definition, the value of this functor on $\B G$ inside $\Glo_{\ab\Lie}$ is precisely the internal mapping object of $\PreAb(\s)$
\[\Ell_{\FF/\s}^{\B G} = \MAP_{\PreAb(\s)}(\underline{\widehat{\B G}},\FF).\]
Let us write $\underline{\widehat{\B G}} = \underline{\widehat{\B G}}_\s$ to emphasise its dependence on $\s$. To see that this assignment defines a functor between Cartesian fibrations over $\Stk$, we straighten these fibrations, which leaves us to define a natural transformation
\[\PreAb\Longrightarrow\Fun(\Glo_{\ab\Lie}, \Shv(-)) \colon \Stk^\op\to \Cat.\]
In other words, we want to check that the internal mapping object assignment above is natural in $\Stk$. This is clear, as for each map of stacks $f\colon \s'\to \s$ of almost finite presentation (afp), the diagram
\[\begin{tikzcd}
{\PreAb(\s)}\ar[d, "{f^\ast}"]\ar[r, "{\Ell_{-/\s}}"]	&	{\Fun(\Glo_{\ab\Lie}, \Shv(\s))}\ar[d, "{f^\ast}"]	\\
{\PreAb(\s')}\ar[r, "{\Ell_{-/\s'}}"]				&	{\Fun(\Glo_{\ab\Lie}, \Shv(\s'))}
\end{tikzcd}\]
naturally commutes. Indeed, the internal mapping objects
\[\MAP_{\PreAb(\s')}(\underline{\widehat{\B G}}_{\s'}, f^\ast \FF) \qquad f^\ast \MAP_{\PreAb(\s)}(\underline{\widehat{\B G}}_{\s},\FF)\]
naturally agree by definition of the functor $f^\ast\colon \Shv(\s)\to \Shv(\s')$---as $f$ is of afp, this morphism lives in the big \'{e}tale site for $\s$, so the evaluation of the sheaf $f^\ast \GG$ on an object $\t'\to\s'$ in the big \'{e}tale site for $\s'$ is the value of $\GG$ on $\t'\to\s'\to\s$ in the big \'{e}tale site for $\s$.
\end{proof}

A key feature of $\Ell^\GM$ is its value when restricted to \emph{abelian} compact Lie groups.

\begin{prop}\label{elliptichomcomparison}
There is a natural equivalence between the composition
\[\int_{\Stk}\PreAb\xrightarrow{\Ell^\GM} \int_{\Stk}\Fun^L(\OS_\Lie,\Shv(-))\to \int_{\Stk} \Fun^L(\Glo_{\ab\Lie}, \Shv(-))\]
induced by restriction along the embedding $\Glo_{\ab\Lie}\to \OS_{\Lie}$, and the assignment sending a pair $\FF/\s$ to the functor
\[\B G\mapsto \HOM_\s(\underline{G}^\vee,\FF).\]
\end{prop}

Note that the homomorphism sheaf assignment above defines a functor on $\int_{\Stk}\PreAb$ as such sheaves satisfy base-change by construction.

\begin{proof}
This identification of functors can be found in \cite[Pr.3.15 \& Ex.6.2]{davidandlennart}. To check for naturality, we recount their argument. Using the constant sheaf--evaluation adjunction used in (\ref{functorgivingellipticcohom}), the natural equivalence $\widehat{\B G}\simeq B \T\times G^\vee$ of (\ref{decompositionofdualofbg}), and that the left adjoint to the forgetful functor $\PreAb(\Spc)\to \Ab(\Spc)$ is given by $B \T\times -$, we obtain the desired object-wise conclusion:
\[\MAP_{\PreAb(\s)}(\underline{\widehat{\B G}},\FF) \simeq \MAP_{\PreAb(\Spc)}(B\T \times G^\vee, \FF(\s))\]
\[ \simeq \MAP_{\Ab(\Spc)}(G^\vee, \FF(\s)) \simeq \MAP_{\Ab(\s)}(\underline{G}^\vee, \FF) = \HOM_\s(\underline{G}^\vee,\FF)\]
The left adjoint used above is natural in the stack variable $\Stk$, meaning it is the value of the left adjoint to the forgetful functor $\int_{\Stk}\PreAb\to \int_{\Stk}\Ab$ at the stack $\s$. Indeed, there is a canonical functor $\int_{\Stk} \PreAb \to \int_{\Stk} \Ab$ associated with the natural transformation $\PreAb(-) \Rightarrow \Ab(-) \colon \Stk^\op \to \Cat$. As the left adjoint above exists on each fibre, and the canonical functor above preserves Cartesian morphisms over $\Stk$, the desired left adjoint above exists by \cite[Pr.7.3.2.6]{ha}. This finishes the proof.
\end{proof}

We are now ready to define the functor, a slight variation of Gepner--Meier's $\Ell^{\GM}$, to be used to compare with Lurie's tempered cohomology.

\begin{mydef}\label{definitionofellgeometric}
Write $\Ell$ for the following variation on $\Ell^{\GM}$, given as the composition
\[\Ell\colon \int_{\Stk}\PreAb \xrightarrow{\Ell^\GM|_{\Glo_\ab}} \int_{\Stk}\Fun(\Glo_\ab,\Shv(-)) \xrightarrow{\mathrm{Lan}}  \int_{\Stk} \Fun^L(\OS,\Shv(-))\]
by restricting $\Ell^\GM$ to $\Glo_\ab$ in the codomain and then taking left Kan extensions along the inclusions $\Glo_\ab\to \Glo \to \OS$. We will also use the notation $\Ell$ for the restriction of this functor to $\int_{\Stk} \PregeoAb$, where $\PregeoAb(\s)$ is the subcategory of preoriented abelian sheaves on $\s$ which are also geometric; this makes sense as the property of being geometric is closed under base-change, which in turn follows as this is true for $\P$-divisible groups.
\end{mydef}

It is formal that this construction agrees with the restriction of $\Ell^\GM$ to the category of orbispaces based on finite groups.

\begin{prop}\label{restrictionisrestriction}
Write $i\colon \Glo\to \Glo_{\Lie}$ for the inclusion and $i^\ast\colon \OS\to \OS_{\Lie}$ for the restriction along $i$. Then we have a natural equivalence of functors
\[i^\ast \circ \Ell^{\GM}  \simeq \Ell\colon \int_{\Stk}\PreAb\to \int_{\Stk}\Fun^L(\OS,\Shv(-)).\]
\end{prop}

\begin{proof}
Consider the diagram of categories
\[\begin{tikzcd}
	{\Glo_{\ab}} & \Glo \\
	{\Glo_{\ab\Lie}} & {\Glo_{\Lie}}
	\arrow["g", from=1-1, to=1-2]
	\arrow["h"', from=1-1, to=2-1]
	\arrow["i"', from=1-2, to=2-2]
	\arrow["j", from=2-1, to=2-2]
\end{tikzcd}\]
where all of the above maps are the obvious inclusions of subcategories. Let us fix a stack $\s$ and a preoriented abelian sheaf $\FF$ and write $F$ for the restriction of $\Ell^\GM_{\FF/\s}$ to $\Glo_{\ab\Lie}$. First, notice that $i^\ast j_! F$, the restriction to $\Glo$ of the left Kan extension along $j$ of $F$, agrees with $g_! h^\ast F$, the left Kan extension along $g$ of the restriction of $F$ to $\Glo_\ab$. Indeed, there is a natural transformation $g_!h^\ast F \to i^\ast j_! F$ courtesy of the universal property of left Kan extensions. To see that this is an equivalence, consider the colimit formula for the value of $j_! F(\B G)$ for a finite group $G$
\[j_! F(\B G) = \underset{\B H \to \B G}{\colim} F(\B H)\]
indexed by maps in $\Glo_{\Lie}$ of the form $\B H \to \B G$ where $H$ is an abelian compact Lie group. Such maps uniquely factor through $\B H \to \B\pi_0 H$ as $G$ is discrete, hence through the subcategory $\Glo \to \Glo_{\Lie}$. In particular, the colimit formula form $i^\ast j_! F(\B G)$ matches that of $g_! h^\ast F$.\\

Unravelling, this shows that the restriction of $\Ell^\GM_{\FF/\s}\colon \OS_{\Lie} \to \Shv(\s)$ to $\Glo$ is the left Kan extension of its restriction to $\Glo_\ab$, in particular, it naturally agrees with $\Ell$ restricted to $\Glo$. As colimit preserving functors out of $\OS$ are determined by their restriction to $\Glo$, we see that the restriction of $\Ell^\GM$ to $\OS$ can be naturally identified with $\Ell$.
\end{proof}

To define our algebraic version of equivariant elliptic cohomology, we will show that the value of $\Ell_{\FF/\s}$ on each object of $\Glo_\ab$ is representable by a stack when $\FF$ is geometric.

\begin{prop}\label{restrictionfactorsthroughstackselliptic}
The composite
\[\int_{\Stk}\PregeoAb\xrightarrow{\Ell} \int_{\Stk} \Fun^L(\OS,\Shv(-)) \to \int_{\Stk} \Fun(\Glo_\ab,\Shv(-))\]
of $\Ell$ with the functor induced by restriction along $\Glo_\ab\to \Glo \to \OS$, factors through $\int_{\Stk} \Fun(\Glo_\ab,\Stk^{\afp,\flat}_{/-})$.
\end{prop}

\begin{remark}\label{finitenessforabeliancptlietoo}
If $\x$ is a preoriented abelian stack over a stack $\s$, where we also assume $\x$ is almost of finite presentation and flat over $\s$, then restricting $\Ell^\GM_{\x/\s}$ to $\Glo_{\ab\Lie}$ also factors through $\Stk^{\afp,\flat}_{/\s}$. Indeed, as an abelian compact Lie group $G$ has the form $\T^r \times H$ for some $r\geq 0$ and a finite abelian group $H$, we obtain the identifications of $\s$-sheaves
\[\HOM_\s(\underline{G}^\vee, \x) \simeq \HOM_\s(\underline{\T^r}^\vee, \x)\times_\s \HOM_\s(\underline{H}^\vee, \x) \simeq \x^{\times_\s r}\times_\s  \HOM_\s(\underline{H}^\vee, \xx) .\]
As colimits in $\OS_{\ab\Lie}$ are sent to colimits in $\Shv(\s)$, it further holds that restricting $\Ell^\GM_{\x/\s}$ to the subcategory of $\OS_{\ab\Lie}$ generated by finite products of representables also factors through $\Stk^{\afp,\flat}_{/\s}$.
\end{remark}

\begin{proof}[Proof of \Cref{restrictionfactorsthroughstackselliptic}]
Courtesy of \Cref{elliptichomcomparison,restrictionisrestriction}, the value of $\Ell_{\FF/\s}$ on $\B G$ for $G$ a finite abelian Lie group is represented by the stack $\HOM_\s(\underline{G}^\vee,\FF)$ which is of afp and flat over $\s$ as $\FF$ is geometric. Using \Cref{factoringhomstackslemma}, we see that the homomorphism sheaf above is represented by an $\s$-stack which is flat and almost of finite presentation, it is even finite, hence we are done by \Cref{equivalencewithconnectiveguys}.
\end{proof}

We finish this section with an algebraic variation of equivariant elliptic cohomology.

\begin{mydef}\label{ellipticglobalsectionsdef}
We define $\O_{\Ell}$ by composing the factorisation of \Cref{restrictionfactorsthroughstackselliptic} together with the functor taking underlying morphisms of $\E_\infty$-algebras in $\QCoh(\s)$, followed by a right Kan extension along $\Glo_{\ab}\to \OS_{\ab}$ and then along $\OS_\ab \to \OS$
\[\int_{\Stk} \PreAb\to \int_{\Stk}\Fun(\Glo_{\ab},\Stk^{\afp,\flat}_{/-}) \to \int_{\Stk} \Fun(\Glo_{\ab}^\op,\CAlg(-))\]
\[\xrightarrow{\mathrm{Ran}} \int_{\Stk} \Fun^R(\OS^\op,\CAlg(-)),\]
where again we denote $\CAlg(\QCoh(\s))$ by $\CAlg(\s)$.
\end{mydef}

It is also reasonable to consider an extension of the above definition to compact Lie groups using \Cref{finitenessforabeliancptlietoo}, but this is irrelevant to our comparison.

%%%%%%%%%%%%%%%%%%%%%%%%%%%%%%%%%%%%%%%%%%%%%%%%%%%%%%%%%%%
%%%%%%%%%%%%%%%%%%%%%%%%%%%%%%%%%%%%%%%%%%%%%%%%%%%%%%%%%%%
%%%%%%%%%%%%%%%%%%%%%%%%%%%%%%%%%%%%%%%%%%%%%%%%%%%%%%%%%%%

\section{Functoriality of tempered cohomology}\label{tempcohomsection}

Lurie's construction of equivariant cohomology theories in \cite{ec3} requires the input of a \emph{preoriented $\P$-divisible group}. Recall that in \cite[Con.4.0.3]{ec3}, Lurie constructs a functor
\[ \A_\G\colon \OS_\ab^\op\to \CAlg_\A \]
for each preoriented $\P$-divisible group (\cite[Df.2.6.1]{ec3}) over an $\E_\infty$-ring $\A$. This is done by first constructing a functor $\Ab_\fin\to \CAlg_\A$ from the $1$-category of finite abelian groups by sending each $H$ to the finite flat $\E_\infty$-$\A$-algebra $\O_{\G[H]}$ representing the functor
\[\CAlg_\A\to \Spc\qquad \mathsf{B}\mapsto \Map_{\Mod_\Z}(H,\G(\mathsf{B}));\]
see the first two paragraphs of \cite[\textsection3.5]{ec3}. Lurie shows in \cite[Th.3.5.5]{ec3} that the data of lifting to the functor $\Ab_\fin\to \CAlg_\A$ through the functor
\[\Ab_\fin\to \Glo_\ab^\op\qquad H\mapsto \B H^\vee\]
is precisely the data of a \emph{preorientation} for $\G$ as a $\P$-divisible group.

\begin{mydef}\label{defofpreorientation}
A \emph{preorientation} for a $\P$-divisible group $\G$ over an $\E_\infty$-ring $\A$ is a morphism of the form $B(\Q/\Z)\to \G(\A)$ in $\Tors(\Spc) \subseteq \Ab(\Spc)$; here we view $B(\Q/\Z)$ as an abelian group object in spaces by identifying it with the connective $\Z$-module $\Q/\Z[1]$; see \Cref{exmapleoftorsion}. Write $\PreBT(\A)$ for the category of \emph{preoriented $\P$-divisible groups over $\A$} given by the fibre product
\[\PreBT(\A) = \BT(\A) \underset{\Tors(\Spc)}{\times} \Tors(\Spc)_{B(\Q/\Z)/}.\]
The functoriality of $\BT(\A)$ in $\A$ induces functoriality of $\PreBT(\A)$ in $\A$, which yields a functor $\PreBT(-) \colon \CAlg \to \Cat$, with associated coCartesian fibration $\int^{\CAlg} \PreBT$.
\end{mydef}

Given a preoriented $\P$-divisible group $\G$ over an $\E_\infty$-ring $\A$, we then define $\A_\G$ as the unique continuous extension of the given functor $\Glo_\ab^\op\to \CAlg_\A$ to $\OS_\ab^\op$, in other words, the right Kan extension. The value of $\A_\G$ on $\B H$ for some finite abelian $H$ is given by the $\E_\infty$-$\A$-algebra $\O_{\G[H^\vee]}$.

\begin{remark}
	Preorientations for $\P$-divisible groups and abelian sheaves (à la \Cref{generaldefintiionofpreorirnation}) are different. In particular, the former uses the abelian group object in spaces $B(\Q/\Z)$ and the latter uses $B\T$. There is a comparison though. For a fixed geometric abelian sheaf $\FF$ over a stack $\s$, the fibre sequence
	\[B\Q \to B(\Q/\Z) \to B\T = B^2 \Z \to B^2 \Q,\]
	induces a fibre sequence of mapping spaces
	\[\Map_{\Ab(\Spc)}(B^2 \Q, \FF(\s)) \to \Map_{\Ab(\Spc)}(B^2 \Z, \FF(\s)) \xrightarrow{\rho} \]
	\[\Map_{\Ab(\Spc)}(B(\Q/ \Z), \FF(\s)) \simeq \Map_{\Ab(\Spc)}(B(\Q/ \Z), \FF[\P^\infty](\s)) \to \Map_{\Ab(\Spc)}(B \Q, \FF(\s)),\]
	the equivalence using the adjunction of \Cref{naturalityoftorsionfunctor} as $B(\Q/\Z)$ is torsion. In particular, the obstructions to lifting a preorientation for $\FF[\P^\infty]$ to a preorientation for $\FF$, as well as the uniqueness of such a lift, are encoded in some rational datum concerning $\FF$; the upper-left and lower-right mapping spaces above. We will see in \Cref{equivalenceofpreorientations} that the torsion object associated with $B^2 \Z = B \T$ is $B(\Q/\Z)$; applying the functor $[\P^\infty](-)$ is another way to define the map $\rho$ in the fibre sequence above.
\end{remark}

Write $\CAlg_{-} \colon \CAlg \to \Cat$ for the functor sending an $\E_\infty$-ring $\A$ to $\CAlg_{\A}$ and a morphism of $\E_\infty$-rings $f\colon \A \to \A'$ to the base change functor $f^\ast \colon \CAlg_{\A} \to \CAlg_{\A'}$ defined by $f^\ast(\mathsf{B}) = \mathsf{B} \otimes_\A \A'$; the left adjoint to the forgetful functor $f_\ast$.

\begin{prop}\label{functorialtemp}
There is a functor
\[\Ga\Temp^{\aff}\colon \int^{\CAlg}\PreBT\to \int^{\CAlg} \Fun^R(\OS_\ab^\op,\CAlg_{-})\]
whose value on a pair $\G/\A$ is $\A_\G$ of Lurie.
\end{prop}

\begin{proof}
As in the proof of \Cref{functorialell}, by taking right Kan extensions, we are left to construct a functor
\[\Ga\Temp^\aff|_{\Glo_\ab}\colon \int^{\CAlg} \PreBT \to \int^{\CAlg} \Fun(\Glo_\ab^\op,\CAlg_{-}).\]

By straightening out these coCartesian fibrations, we are left to show that Lurie's tempered cohomology construction yields a natural transformation
\[\PreBT\Longrightarrow \Fun(\Glo_\ab^\op,\CAlg_{-}) \colon \CAlg\to \Cat.\]

To this end, we note that for a morphism of $\E_\infty$-rings $f\colon \A \to \A'$, the diagram of categories
\[\begin{tikzcd}
{\PreBT(\A)}\ar[d, "{f^\ast}"]\ar[r, "{\A_{-}}"]	&	{\Fun(\Glo_\ab^\op, \CAlg_{\A})}\ar[d, "{f^\ast}"]	\\
{\PreBT(\A')}\ar[r, "{\A_{-}}"]			&	{\Fun(\Glo_\ab^\op, \CAlg_{\A'})}
\end{tikzcd}\]
naturally commutes, as there is a natural equivalence of $\E_\infty$-$\A'$-algebras
\[f^\ast  \O_{\G[H^\vee]} = \A'\otimes_{\A} \O_{\G[H^\vee]} \simeq \O_{f^\ast \G[H^\vee]}\]
for every finite abelian group $H$, by the corepresentability and universal properties of both sides.
\end{proof}

Our goal is now to extend Lurie's tempered cohomology construction from affine stacks to general stacks. We begin with preorientations.

\begin{mydef}\label{preorientationoverstack}
 A \emph{preorientation} of a $\P$-divisible group $\G$ over a stack $\s$ is then a map in $\Tors(\Spc)$ of the form $B(\Q/\Z)\to \G(\s)$. As $B(\Q/\Z)$ lies in $\Mod_\Z^{\cn,\Tors}$, we see that the space of preorientations for $\G$ over $\s$ can be written as the limit
\[\Pre(\G)\simeq \lim_{f\colon \Spec \A\to \s \in \s_\et^\aff} \Pre(f^\ast \G)\]
using (\ref{globalsectionsofpdiv}). Let us write $\PreBT(\s) = \lim_{\s_\et^\aff} \PreBT(\A)$ for the category of preoriented $\P$-divisible groups.
\end{mydef}

In particular, notice that if $\s=\Spec \A$ is an affine stack, then a $\P$-divisible group $\G$ over $\s$ matches the definition of Lurie \cite[Df.2.6.1]{ec3} by the Yoneda lemma. Similarly, a preorientation of $\G$ over $\s$ also agrees with the definition of Lurie \cite[Df.2.6.8]{ec3} when $\s$ is affine, as the category of \'{e}tale morphisms from an affine stack into $\s$ now has a terminal object.

\begin{mydef}\label{definitionoftempglobal}
Let $\s$ be a stack. By construction, we can write the symmetric monoidal category $\QCoh(\s)$ of quasi-coherent sheaves on $\s$ as the limit
\[\QCoh(\s)\xrightarrow{\simeq} \lim_{\s_\et^\aff} \QCoh(\Spec \A)\simeq \lim_{\s_\et^\aff} \Mod_\A,\]
the second equivalence follows from the natural symmetric monoidal equivalence between $\QCoh(\Spec \A)$ and $\Mod_\A$; see \cite[Pr.2.2.3.3]{sag}. Fix a preoriented $\P$-divisible group $\G$ over $\s$. For each $f\colon \Spec\A\to \s$ inside $\s_\et^\aff$ we obtain a functor
\[\A_{f^\ast\G}\colon \Glo_\ab^\op\to \CAlg_\A.\]
By \Cref{functorialtemp}, this defines an object inside the limit of functor categories
\[ \Fun(\Glo_\ab^\op,\CAlg(\s))\simeq \Fun(\Glo_\ab^\op,\lim_{\Spec \A\to \s} \CAlg_{\A}) \simeq \lim_{\Spec \A\to \s} \Fun(\Glo_\ab^\op,\CAlg_{\A})\]
where the limits are again indexed over $\s_\et^\aff$. Taking a right Kan extension of this functor defined by this collection of $\{\A_{f^\ast \G}\}$ along the inclusions $\Glo_\ab\to \OS_\ab \to \OS$ yields our desired \emph{algebraic tempered cohomology} functor
\[\O_{\Temp_{\G/\s}}\colon \OS^\op\to \CAlg(\s).\]
Using \Cref{functorialtemp} again, this construction sending a pair $\G/\s$ to the functor $\O_{\Temp_{\G/\s}}$ refines further to a functor
\[\O_{\Temp}\colon \int_{\Stk}\PreBT\to \int_{\Stk} \Fun^R(\OS^\op,\CAlg(-)).\]
\end{mydef}

From the same comments after \Cref{preorientationoverstack}, if $\s=\Spec \A$ is affine, then the global sections of $\O_{\Temp_{\G/\s}}$ evaluated on an orbispace $X$ in $\OS_\ab$ agrees with Lurie's $\A_\G^X$.\footnote{Much of these definitions also work if $\s$ has a smooth or even fpqc-hypercover by affines, as $\M_{\BT}$ is an fpqc-hypersheaf. We will not pursue that direction here; this will appear in future work with Balderrama and Linskens.}\\

Minor changes lead us to our geometric tempered cohomology functor.

\begin{mydef}\label{geometricdefinitionoftemp}
Let $\s$ be a stack and $\G$ be a $\P$-divisible group over $\s$. Notice that for each morphism $f\colon \Spec \A \to \s$ in $\s_\et^\aff$, the restriction of $\A_{f^\ast \G}\colon \OS_\ab^\op \to \CAlg_\A$ to $\Glo_\ab$ factors through $\CAlg_{\A}^{\afp,\flat}$, the category of $\E_\infty$-$\A$-algebras which are flat and of almost finite presentation. Indeed, its value on $\B H$ is the finite flat $\E_\infty$-$\A$-algebra $\O_{f^\ast \G[H^\vee]}$ of \Cref{pdivdef} by \Cref{functorialtemp} and Lurie's construction \cite[Con.4.0.3]{ec3}. Consider the functor 
\begin{equation}\label{temponeachguyinetalesite} \Glo_\ab\xrightarrow{\A_{f^\ast \G}^\op} (\CAlg_{\A}^{\afp,\flat})^{\op}\simeq \Aff^{\afp,\flat}_{/\Spec \A}\to \Shv(\A)\end{equation}
defined by this factorisation and the Yoneda embedding. As in \Cref{definitionoftempglobal}, as $f\colon \Spec \A\to \s$ varies over $\s_\et^\aff$, the assignment sending $f\colon \Spec \A \to \s$ in $\s^\aff_\et$ to the functor (\ref{temponeachguyinetalesite}) defines an object in the functor category
\[\Fun(\Glo_\ab,\Shv(\s))\simeq \Fun(\Glo_\ab,\lim_{\s_\et^\aff}\Shv(\A)) \simeq \lim_{\s_\et^\aff} \Fun(\Glo_\ab,\Shv(\A)).\]
Moreover, by \Cref{functorialtemp}, the process of taking $\s$ to this functor $\Glo_\ab \to \Shv(\s)$ itself defines a functor
\[\int_{\Stk} \PreBT\to \int_{\Stk} \Fun(\Glo_\ab,\Shv(-))\]
using the facts that the assignment $\s\mapsto \s_\et^\aff$ is functorial in $\s$, as is taking limits. By first left Kan extending from $\Glo_\ab$ to $\Glo$, and then left Kan extending again along $\Glo\to \OS$, we obtain our \emph{geometric tempered cohomology} functor
\[\Temp\colon \int_{\Stk} \PreBT \to \int_{\Stk} \Fun^L(\OS,\Shv(-)).\]
\end{mydef}

By construction, the value of $\Temp_{\G/\s}$ on $\B H$ is the homomorphism sheaf of \Cref{defhomsheafforpdiv}.

\begin{prop}\label{comparisontohomstacktempered}
There is a natural equivalence between the composition
\[\int_{\Stk}\PreBT\xrightarrow{\Temp} \int_{\Stk}\Fun^L(\OS,\Shv(-))\to \int_{\Stk} \Fun(\Glo_\ab, \Shv(-))\]
induced by the embedding $\Glo_\ab\to \OS$, and the assignment sending a pair $\G/\s$ to the functor
\[\B H\mapsto \HOM_\s(\underline{H}^\vee,\G).\]
\end{prop}

\begin{proof}
The natural equivalence $\Temp_{\G/\s}(\B G) \simeq \HOM_\s(\underline{H}^\vee,\G)$ of $\s$-sheaves will be given by taking a limit of natural equivalences $\Temp_{f^\ast\G/\A}(\B G) \simeq \HOM_\A(\underline{H}^\vee,f^\ast\G)$ for all étale morphisms $\Spec \A \to \s$, so we are reduced to the affine case. The conclusion then follows, as by \Cref{geometricdefinitionoftemp}, $\Temp_{\G/\Spec \A}(\B H)$ is defined on affines by $\O_{\G[H^\vee]}$, which corepresents $\HOM_\A(\underline{H}^\vee, \G)$ by \Cref{defhomsheafforpdiv}.
\end{proof}

%%%%%%%%%%%%%%%%%%%%%%%%%%%%%%%%%%%%%%%%%%%%%%%%%%%%%%%%%%%
%%%%%%%%%%%%%%%%%%%%%%%%%%%%%%%%%%%%%%%%%%%%%%%%%%%%%%%%%%%
%%%%%%%%%%%%%%%%%%%%%%%%%%%%%%%%%%%%%%%%%%%%%%%%%%%%%%%%%%%

\section{Proving the comparison theorems}\label{comparisonsection}
We almost have enough mathematics to prove \Cref{maintheoremalg,maintheoremgeo}. First, we need to compare a vital piece of data used to define these theories: \emph{preorientations}.

\begin{prop}\label{equivalenceofpreorientations}
There is an equivalence of abelian group objects in spaces
\[B \T[\P^\infty] \simeq B (\Q/\Z).\]
In particular, there is a functor
\[[\P^\infty] \colon \int_{\Stk} \PregeoAb \to \int_{\Stk}\PreBT\]
refining that of \Cref{firstattemptatfunctorusingpinfty}.\footnote{For an abelian variety $\x$, Lurie defines another construction of $\x[\P^\infty]$ for $\s=\Spec \A$ and $\A$ connective, simply as the functor $\CAlg_{\A}^\cn\to \Mod^\cn_\Z$ sending $\b$ to the fibre of $\x[\P^\infty](\b)\to \x[\P^\infty](\b)\otimes \Q$; see \cite[Con.2.9.1]{ec3}. These two definitions match when $\s=\Spec\A$ is affine, as in both cases $\x[\P^\infty](\b)$ is the terminal torsion connective $\Z$-module over $\x(\b)$.} 
\end{prop}

\begin{proof}
The equivalence in abelian group objects in spaces comes from \cite[Rmk.6.4.8]{ec1}. Indeed, first, we note that both $B\T[\P^\infty]$ and $B(\Q/\Z)$ are torsion objects, meaning they are in the essential image of the left adjoint $\Tors(\Spc)\to \Ab(\Spc)$ to the functor $[\P^\infty]$ of \Cref{naturalityoftorsionfunctor}. Phrased differently, under the equivalences of categories
\[\Tors(\Spc)\simeq \Mod_\Z^{\cn,\Tors}\to \Mod_\Z^\cn\simeq \Ab(\Spc)\]
of \cite[Ex.3.5.10]{ec3}, the connective $\Z$-module spectra associated with these objects of $\Ab(\Spc)$ are both torsion. To see these two objects agree as $\P$-torsion objects, it now suffices to check they agree when evaluated on cyclic groups. By \cite[Rmk.6.4.8]{ec1}, the value of $B\T[\P^\infty]$ on the finite abelian group $\Z/n\Z$ for some $n\geq 1$, can be calculated as the fibre of the $n$-fold multiplication map $[n]\colon B\T(\Z)\to B\T(\Z)$, also known as the \emph{$n$-torsion of $B \T(\Z)$}; see \Cref{def:ntorsion}. The values of $B\T[\P^\infty]$ and $B(\Q/\Z)$ on these cyclic groups are now clearly both the discrete $\Z$-module $B (\Z/n\Z)$. To define the desired functor, we need to define a natural transformation of the functors
\[\PregeoAb(-) \Rightarrow \PreBT(-)\colon \Stk^\op \to \Cat.\]
As in the proof of \Cref{firstattemptatfunctorusingpinfty}, for each stack $\s$, we have $\PreBT(\s) = \lim_{\s_\et^\aff} \PreBT(\A)$ by definition, so it suffices to construct the desired natural transformation when restricting these functors to affine stacks. In this case, we note that the functoriality of $[\P^\infty]$ courtesy of \Cref{naturalityoftorsionfunctor} gives us a functor
\[\PreAb(\Spec \A) \simeq \Ab(\Spec \A)_{\underline{B\T}/} \xrightarrow{[\P^\infty]_{\Shv(\Spec \A)}} \Tors(\Spec \A)_{\underline{B(\Q/\Z)}/} \]
natural in $\Spec \A$. Restricting to $\PregeoAb(\Spec \A)$ factors through $\BT(\A)_{\underline{B(\Q/\Z)}/} \simeq \PreBT(\A)$, yielding our desired natural transformation.
\end{proof}

We can finally define the last of our functors needed for our comparison.

\begin{notation}
Define the functors $\O_{\Temp_{[\P^\infty]}}$ and $\Temp_{[\P^\infty]}$ as the following composites:
\[\O_{\Temp_{[\P^\infty]}} \colon \int_{\Stk} \PregeoAb\xrightarrow{[\P^\infty]} \int_{\Stk} \PreBT \xrightarrow{\O_{\Temp}} \int_{\Stk} \Fun^R(\OS^\op,\CAlg(-))\]
\[\Temp_{[\P^\infty]} \colon \int_{\Stk} \PregeoAb\xrightarrow{[\P^\infty]} \int_{\Stk} \PreBT \xrightarrow{\Temp} \int_{\Stk} \Fun^L(\OS,\Shv(-))\]
\end{notation}

To prove \Cref{maintheoremalg,maintheoremgeo}, we start with the crucial special case.

\begin{prop}\label{specialcase}
There is a natural equivalence
\[\Temp_{[\P^\infty]}|_{\Glo_\ab} \xRightarrow{\simeq} \Ell|_{\Glo_\ab} \colon \int_{\Stk} \PregeoAb\to \int_{\Stk} \Fun(\Glo_\ab,\Shv(-))\]
induced by the natural map $\FF[\P^\infty]\to \FF$.
\end{prop}

\begin{proof}
By \Cref{elliptichomcomparison} and \Cref{comparisontohomstacktempered}, we see that the two above functors are equivalent to the assignments sending a pair $\FF/\s$ to the functors sending a finite abelian group $G$ to the homomorphism stacks 
\[\Ell_{\FF/\s}^{\B G}\simeq\HOM_\s(\underline{G}^\vee,\FF)\qquad \Temp_{\FF[\P^\infty]/\s}^{\B G}\simeq\HOM_\s(\underline{G}^\vee,\FF[\P^\infty])\]
respectively. The natural transformation between these functors is induced by the natural map $\varepsilon \colon \FF[\P^\infty]\to \FF$, which is itself natural in $G$, and is also natural in $\int_{\Stk}\PregeoAb$ as this map is a counit of the natural adjunction from \Cref{naturalityoftorsionfunctor}. By \Cref{factoringhomstackslemma}, the map $\varepsilon$ induces our desired natural equivalence between these homomorphism functors:
\[\HOM_\s(\underline{G}^\vee,\FF[\P^\infty])\xrightarrow{\simeq} \HOM_\s(\underline{G}^\vee,\FF) \qedhere \]
\end{proof}

\Cref{maintheoremalg,maintheoremgeo} now follow.

\begin{proof}[Proof of \Cref{maintheoremgeo}]
Both of the functors $\Ell$ and $\Temp_{[\P^\infty]}$ are left Kan extended from their restrictions along the inclusions $\Glo_\ab\to \Glo\to \OS$. We are then left to show that these restrictions to $\Glo_\ab^\op$ are naturally equivalent, which is precisely \Cref{specialcase}.
\end{proof}

\begin{proof}[Proof of \Cref{maintheoremalg}]
By definition, both $\O_{\Ell}$ and $\O_{\Temp_{[\P^\infty]}}$ are defined by right Kan extending their restrictions to $\Glo_\ab^\op$. On the other hand, by \Cref{elliptichomcomparison} and \Cref{specialcase}, the restrictions of $\Ell$ and $\Temp$ to $\Glo_\ab$ factor through $\int_{\Stk} \Fun(\Glo_\ab,\Stk^{\afp,\flat}_{/-})$, and so we are done as the associated algebras of quasi-coherent sheaves of these restrictions are precisely the restrictions of $\O_{\Ell}$ and $\O_{\Temp_{[\P^\infty]}}$.
\end{proof}

%%%%%%%%%%%%%%%%%%%%%%%%%%%%%%%%%%%%%%%%%%%%%%%%%%%%%%%%%%%
%%%%%%%%%%%%%%%%%%%%%%%%%%%%%%%%%%%%%%%%%%%%%%%%%%%%%%%%%%%
%%%%%%%%%%%%%%%%%%%%%%%%%%%%%%%%%%%%%%%%%%%%%%%%%%%%%%%%%%%

\section{An application to equivariant topological modular forms}\label{dualisableiltysection}

As an example application of \Cref{maintheoremalg,maintheoremgeo}, let us repeat and refine the argument outlined in \cite[Rmk.10.6]{davidandlennart}, to show that $\TMF^{\B G}$ is a dualisable $\TMF$-module for certain compact Lie groups $G$; in \emph{ibid}, only a plausibility argument could be given as no explicit comparison was constructed. Recall the definitions of \emph{orientations} from \cite[Df.5.17]{davidandlennart} and \cite[\textsection4]{ec2}. We will also write the functor $\Ell^\GM$ of \Cref{functorialell} as $\Ell$ to declutter some notation.

\begin{theorem}[{\Cref{dualisableproperty}}]\label{dualisableforfinitegroups}
Let $G$ be a compact Lie group which can be written as the product of a torus with a finite group and $\x$ be an \textbf{oriented} elliptic curve over a stack $\s$. Then the $\E_\infty$-$\O_\s$-algebra $\O_{\Ell_{\x/\s}^{\B G}}$ is dualisable as a quasi-coherent sheaf over $\s$. In particular, $\TMF^{\B G}$ is dualisable as a $\TMF$-module.
\end{theorem}

\begin{proof}
The ``in particular'' case follows by setting $\s$ to be the moduli stack of oriented elliptic curves $\M_\Ell^\ori$ and $\x$ to be the universal oriented elliptic curve $\e^\ori$ (\cite[Ex.10.4]{davidandlennart} and \cite[\textsection7]{ec2}) and the fact that the global sections functor, which sends the structure sheaf of ${\Ell_{\x/\M_\Ell^\ori}^{\B G}}$ to $\TMF^{\B G}$ in this case, induces a symmetric monoidal equivalence $\QCoh(\M_\Ell^{\mathrm{or}})\simeq \Mod_{\TMF}$; see \cite[Th.7.2]{akhilandlennart} or \cite[Th.A]{reconstruction}.\footnote{This argument also shows that, given the hypotheses of \Cref{dualisableforfinitegroups}, $\Ga\Ell_{\x/\s}^{\B G}$ is dualisable as a $\Ga(\s)$-module if the defining map $\s\to \M_\Ell^\ori$ of $\x$ is quasi-affine, or at least such that the global sections functor $\Ga\colon \QCoh(\s) \to \Mod_{\Ga(\s)}$ is an equivalence; see \cite{reconstruction} for more examples of such stacks.}\\

In general, by definition as a left Kan extension, the value of ${\Ell^{\B G}_{\x/\s}}$ is given by the colimit of ${\Ell^{\B H}_{\x/\s}}$ over all maps $\B H\to \B G$ in $\Glo_{\Lie}$, where $G$ is our given compact Lie group and $H$ is an abelian compact Lie group. To simplify this indexing category, note that as there is a splitting $G \simeq G' \times\T^r$ for some finite group $G'$ and some nonnegative integer $r\geq 0$, we have a chosen equivalence $\B G \simeq \B G' \times \B\T^r $ inside $\Glo_{\Lie}$. The indexing category in question then splits into a product of the indexing categories of maps $\B H' \to \B G'$ and $\B H'' \to \T^r$ in $\Glo_{\Lie}$ with both $H'$ and $H''$ being abelian compact Lie groups. The maps $\B H' \to \B G'$ in the first of these categories uniquely factor through $\B H' \to \B \pi_0 H'$ as $G'$ is a finite group, as in the proof of \Cref{restrictionisrestriction}. In particular, this first category is equivalent to that of maps $\B H' \to \B G'$ in $\Glo$, where $H'$ is now a finite abelian group. The second of these categories has a final object $\B H'' = \B\T^r$. We now have the equivalences of abelian stacks
\[ \Ell^{\B G}_{\x/\s} \simeq \colim_{\B H \to \B G} \Ell^{\B H}_{\x/\s} \simeq \colim_{\B H' \to \B G'} \Ell^{\B (H'\times \T^r)}_{\x/\s}.\]
From the definition of $\Ell^{\B(H'\times \T^r)}_{\x/\s}$ as a homomorphism stack, the base-change property for such stacks, and the equivalence $\Ell_{\x/\s}^{\B \T^r} \simeq \x^{\times_\s r}$ of \Cref{finitenessforabeliancptlietoo}, we can compute
\[\Ell^{\B(H'\times \T^r)}_{\x/\s} \simeq \HOM_\s(\underline{(H'\times \T^r)}^\vee, \x) \simeq \HOM_\s(\underline{H'}^\vee, \x) \times_\s \HOM_\s(\underline{\Z^r}, \x)\]
\[ \simeq \HOM_\s(\underline{H'}^\vee, \x) \times_\s \x^{\times_\s r} \simeq \HOM_{\x^{\times_\s r}}(\underline{H'}^\vee, \x^{\times_\s (r+1)}) \simeq \Ell^{\B H'}_{\x^{\times_\s (r+1)}/\x^{\times_\s r}}.\]
Combining these chains of equivalences together, we obtain equivalences of stacks
\[\Ell^{\B G}_{\x/\s} \simeq \colim \Ell^{\B H'}_{\x^{\times_\s (r+1)}/\x^{\times_\s r}} \simeq \Ell^{\B G'}_{\x^{\times_\s (r+1)}/\x^{\times_\s r}}.\]
Taking structure sheaves, we are reduced to the case where $G$ is a finite group as $\x^{\times_\s (r+1)} \to \x^{\times_\s r}$ is an oriented elliptic curve by base-change. In this case, the equivalence of \Cref{restrictionisrestriction} between $\Ell^\GM$ and $\Ell$ allows us to appeal to \Cref{maintheoremalg}, which yields an equivalence of $\E_\infty$-algebras in $\QCoh(\s)$ of the form $\O_{\Ell^{\B G}_{\x/\s}}\simeq \O_{\Temp^{\B G}_{\x[\P^\infty]/\s}}$. Now take an affine étale hypercover $\Spec \A^\bullet\to \s$. By construction, the above $\E_\infty$-algebra in $\QCoh(\s)$ is the cosimplicial limit
\[\O_{\Temp^{\B G}_{\x[\P^\infty]/\s}} \simeq \lim \O_{\Temp^{\B G}_{\x_{\A^\bullet}[\P^\infty]/\A^\bullet}},\]
where $\x_{\A^\bullet}$ denotes the base-change of $\x$ to $\Spec \A^\bullet$. By \cite[Pr.4.6.1.11]{ha} and the equivalence $\CAlg(\QCoh(\s))\simeq \lim \CAlg_{\A^\bullet}$, it suffices to see that each $\E_\infty$-$\A^\bullet$-algebra given by the global sections of $\O_{\Temp^{\B G}_{\x_{\A^\bullet}[\P^\infty]/\A^\bullet}} = (\A^\bullet)^{\B G}_{\x_{\A^\bullet}[\P^\infty]}$ is dualisable as an $\A^\bullet$-module, hence we are reduced to the affine case $\s=\Spec \A$. This now follows from \cite[Pr.4.7.9]{ec3}, which states that $\A_{\x[\P^\infty]}^{\B G}$ is perfect, hence also dualisable by \cite[Prs.7.2.4.2 \& 7.2.4.4]{ha}, as an $\A$-module.
\end{proof}

\begin{remark}
	In fact, one can show that in the affine case appearing at the end of the above proof, that $\A_{\x[\P^\infty]}^{\B G}$ is actually \emph{self-dual} as an $\A$-modules. Indeed, this is shown in \cite[Rmk.10.6]{davidandlennart}, and we paraphrase their argument here. For each finite group $G$, Lurie defines an $\A$-linear stable category of local systems $\mathrm{LocSys}_{\x[\P^\infty]}(\B G)$ with functorial restrictions $f^\ast\colon \mathrm{LocSys}_{\x[\P^\infty]}(\B H)\to \mathrm{LocSys}_{\x[\P^\infty]}(\B G)$ for each map $f\colon \B G\to \B H$ in $\Glo$. If $G$ is the trivial group $e$, then $\mathrm{LocSys}_{\x[\P^\infty]}(\ast)$ is canonically equivalent to $\Mod_{\A}$ and the functor $p_\ast$ from $\mathrm{LocSys}_{\x[\P^\infty]}(\B G)$ to $\mathrm{LocSys}_{\x[\P^\infty]}(\ast)$ right adjoint to the restriction $p^\ast$, induced by the projection $p\colon \B G\to \ast$, sends the ``constant local system'' $\underline{\A}_{\B G}$ to the tempered cohomology theory $\A_{\x[\P^\infty]}^{\B G}$. Lurie's tempered ambidexterity theorem \cite[Th.1.1.21]{ec3} states the existence of a canonical equivalence between this right adjoint $f_\ast$ and the left adjoint of $f^\ast$, denoted by $f_!$. Using this equivalence, these adjunctions, and the universal property of ``constant local systems'', see \cite[Rmk.5.1.20]{ec3}, we obtain the desired (self-) duality:
\[\MAP_\A(\A_{\x[\P^\infty]}^{\B G},\A)\simeq \MAP_\A(f_\ast \underline{\A}_{\B G},\A) \simeq \MAP_\A(f_! \underline{\A}_{\B G},\A) \]
\[\simeq \MAP_{\mathrm{LocSys}_{\x[\P^\infty]}(\B G)}(\underline{\A}_{\B G},f_\ast \A)\simeq f_\ast \underline{\A}_{\B G}\simeq \A^{\B G}_{\x[\P^\infty]}.\]
\end{remark}

%%%%%%%%%%%%%%%%%%%%%%%%%%%%%%%%%%%%%%%%%%%%%%%%%%%%%%%%%%%
%%%%%%%%%%%%%%%%%%%%%%%%%%%%%%%%%%%%%%%%%%%%%%%%%%%%%%%%%%%
%%%%%%%%%%%%%%%%%%%%%%%%%%%%%%%%%%%%%%%%%%%%%%%%%%%%%%%%%%%

%REFERENCES!!!
\addcontentsline{toc}{section}{References}
\scriptsize
%\bibliography{C:/Users/jackd/Dropbox/Work/references} 

\bibliography{/Users/jackdavies/Dropbox/Work/references} 

\bibliographystyle{alpha}

%\newpage
%\listoftodos

\end{document}